\documentclass[a4paper]{amsart}
\usepackage[T1]{fontenc}
\usepackage[latin1]{inputenc}
\usepackage{amssymb}
\usepackage{stmaryrd}
\usepackage{enumerate}
\usepackage[lite]{amsrefs}

\newcommand*{\abs}[1]{\lvert#1\rvert}
\newcommand*{\norm}[1]{\lVert#1\rVert}

\newcommand*{\prto}{\twoheadrightarrow}

\newcommand*{\biinv}{{\sslash}}

\newcommand*{\defeq}{\mathrel{:=}}

\newcommand*{\brd}{-\hspace{0pt}}
\newcommand*{\nbd}{\nobreakdash-\hspace{0pt}}

\newcommand*{\hot}{\mathbin{\hat{\otimes}}}
\newcommand*{\barotimes}{\mathbin{\bar{\otimes}}}
\newcommand*{\Lhot}{\mathbin{\hat{\otimes}}^{\mathbb{L}}}
\newcommand*{\congto}{\overset{\cong}\to}
\newcommand*{\conj}[1]{\overline{#1}}

\newcommand*{\Mod}{\mathsf{Mod}}
\newcommand*{\Der}{\mathsf{Der}}
\newcommand*{\Extensions}{\mathsf{Ext}}
\newcommand*{\Fine}{\mathsf{Fine}}
\newcommand*{\Vonn}{\mathsf{vN}}
\newcommand*{\COM}{\mathsf{CO}}

\newcommand*{\alg}{\mathrm{alg}}
\newcommand*{\fin}{\mathrm{fin}}
\newcommand*{\ID}{\mathrm{id}}

\newcommand*{\BT}{\mathcal{BT}}
\newcommand*{\Sch}{\mathcal{S}}

\newcommand*{\Hecke}{\mathcal{H}}
\newcommand*{\Controlled}{\mathcal{C}}
\newcommand*{\T}{{\mathcal{T}}}
\newcommand*{\Tc}{{\mathcal{T}}_{\mathrm{c}}}

\DeclareMathOperator{\supp}{supp}
\DeclareMathOperator{\Map}{Map}
\DeclareMathOperator{\stab}{Stab}
\DeclareMathOperator{\Fix}{Fix}
\DeclareMathOperator{\vol}{vol}
\DeclareMathOperator{\Hom}{Hom}
\DeclareMathOperator{\Endo}{End}
\DeclareMathOperator{\Aut}{Aut}
\DeclareMathOperator{\Right}{\mathbb{R}}

\DeclareMathOperator{\Ext}{Ext}
\DeclareMathOperator{\Tor}{Tor}
\DeclareMathOperator{\rk}{rk}
\DeclareMathOperator{\K}{K}
\DeclareMathOperator{\HH}{HH}
\DeclareMathOperator{\Rep}{Rep}
\DeclareMathOperator{\cInd}{\text{c-Ind}}

\DeclareMathOperator{\tr}{tr}
\DeclareMathOperator{\Eul}{Eul}
\DeclareMathOperator{\clos}{cl}
\DeclareMathOperator{\size}{size}

\newcommand*{\EP}{{\mathrm{EP}}}
\newcommand*{\univ}{{\mathrm{univ}}}
\newcommand*{\op}{\mathrm{op}}
\newcommand*{\U}{{\mathrm{U}}}

\newcommand*{\Gss}{{G_{\mathrm{ss}}}}

\newcommand*{\C}{{\mathbb{C}}}
\newcommand*{\R}{{\mathbb{R}}}
\newcommand*{\Z}{{\mathbb{Z}}}
\newcommand*{\N}{{\mathbb{N}}}
\newcommand*{\Q}{{\mathbb{Q}}}

\theoremstyle{plain}
\newtheorem{theorem}{Theorem}
\newtheorem{proposition}[theorem]{Proposition}
\newtheorem{lemma}[theorem]{Lemma}

\theoremstyle{definition}
\newtheorem{definition}[theorem]{Definition}

\begin{document}

\title[Homological algebra for Schwartz algebras]{Homological algebra for
  Schwartz algebras of reductive p-adic groups}
\author{Ralf Meyer}
\email{rameyer@math.uni-muenster.de}

\address{Mathematisches Institut\\
         Westfälische Wilhelms-Universität Münster\\
         Einsteinstr.\ 62\\
         48149 Münster\\
         Germany
}


\thanks{This research was supported by the EU-Network \emph{Quantum
    Spaces and Noncommutative Geometry} (Contract HPRN-CT-2002-00280)
  and the \emph{Deutsche Forschungsgemeinschaft} (SFB 478).}

\begin{abstract}
  Let~\(G\) be a reductive group over a non-Archimedean local field.  Then the
  canonical functor from the derived category of smooth tempered
  representations of~\(G\) to the derived category of all smooth
  representations of~\(G\) is fully faithful.  Here we consider
  representations on bornological vector spaces.  As a consequence, if~\(G\)
  is semi-simple, \(V\) and~\(W\) are tempered irreducible representations
  of~\(G\), and \(V\) or~\(W\) is square-integrable, then
  \(\Ext_G^n(V,W)\cong0\) for all \(n\ge1\).  We use this to prove in full
  generality a formula for the formal dimension of square-integrable
  representations due to Schneider and Stuhler.
\end{abstract}
\maketitle

\section{Introduction}
\label{sec:intro}

Let~\(G\) be a linear algebraic group over a non-Archimedean local field whose
connected component of the identity element is reductive; we briefly call such
groups \emph{reductive \(p\)-adic groups}.  For the purposes of exposition, we
assume throughout the introduction that the connected centre of~\(G\) is
trivial, although we treat groups with arbitrary centre in the main body of
this article.

We are going to compare homological and cohomological computations for the
Hecke algebra \(\Hecke(G)\) and the Harish-Chandra Schwartz algebra
\(\Sch(G)\).  Our main result asserts that the derived category of \(\Sch(G)\)
is a full subcategory of the derived category of \(\Hecke(G)\).  These derived
categories incorporate a certain amount of functional analysis because
\(\Sch(G)\) is more than just an algebra.  Before we discuss this, we sketch
two purely algebraic applications of our main theorem.

Let \(\Mod_\alg(G)\) be the category of smooth representations of~\(G\) on
complex vector spaces.  We compute some extension spaces in this Abelian
category.  If both \(V\) and~\(W\) are irreducible tempered representations
and one of them is square-integrable, then \(\Ext^n_G(V,W)=0\) for \(n\ge1\).
If the local field underlying~\(G\) has characteristic~\(0\), this is proven
by very different means in~\cite{Schneider-Stuhler}.  We get a more
transparent proof that also works in prime characteristic.

The vanishing of \(\Ext^n_G(V,W)\) is almost trivial if \(V\) or~\(W\) is
supercuspidal because then \(V\) or~\(W\) is both projective and injective in
\(\Mod_\alg(G)\).  This is related to the fact that supercuspidal
representations are isolated points in the admissible dual.  Square-integrable
representations are isolated points in the tempered dual.  Hence they are
projective and injective in an appropriate category
\(\Mod\bigl(\Sch(G)\bigr)\) of tempered smooth representations of~\(G\).  Both
\(\Mod\bigl(\Sch(G)\bigr)\) and \(\Mod_\alg(G)\) are full subcategories in a
larger category \(\Mod(G)\).  That is, we have fully faithful functors
\[
\Mod\bigl(\Sch(G)\bigr) \rightarrow \Mod(G) \leftarrow \Mod_\alg(G).
\]
We will show that the induced functors between the derived categories,
\[
\Der\bigl(\Sch(G)\bigr) \rightarrow \Der(G) \leftarrow \Der_\alg(G),
\]
are still fully faithful.  This contains the vanishing result for \(\Ext\) as
a special case.

Another application involves Euler characteristics for square-integrable
representations.  Let~\(V\) be an irreducible square-integrable representation
of~\(G\).  By a theorem of Joseph Bernstein, any finitely generated smooth
representation of~\(G\), such as~\(V\), has a finite type projective
resolution \(P_\bullet\to V\).  Its Euler characteristic is defined as
\[
\chi(V)\defeq \sum (-1)^n [P_n] \in\K_0\bigl(\Hecke(G)\bigr).
\]
Since~\(V\) is square-integrable, it is a projective \(\Sch(G)\)-module and
therefore defines a class \([V]\in\K_0\bigl(\Sch(G)\bigr)\).  We show that the
map \(\K_0\bigl(\Hecke(G)\bigr)\to\K_0\bigl(\Sch(G)\bigr)\) induced by the
embedding \(\Hecke(G)\to\Sch(G)\) maps \(\chi(V)\) to \([V]\).  This is useful
because in~\cite{Schneider-Stuhler} Peter Schneider and Ulrich Stuhler
construct very explicit finite type projective resolutions, so that we get a
nice formula for \([V]\in\K_0\bigl(\Sch(G)\bigr)\).  This implies an explicit
formula for the formal dimension of~\(V\), which is proven
in~\cite{Schneider-Stuhler} if~\(V\) is supercuspidal or if~\(G\) has
characteristic~\(0\).  One consequence of this formula is that the formal
dimensions are quantised, that is, they are all multiples of some
\(\alpha>0\).  This allows to estimate the number of irreducible
square-integrable representations that contain a given representation of a
compact open subgroup of~\(G\).

Although these applications can be stated purely algebraically, their proofs
require functional analysis.  We may view \(\Sch(G)\) just as an algebra and
consider the category \(\Mod_\alg\bigl(\Sch(G)\bigr)\) of modules over
\(\Sch(G)\) in the algebraic sense as in~\cite{Schneider-Zink}.  However, the
functor \(\Der_\alg\bigl(\Sch(G)\bigr)\to\Der_\alg(G)\) fails to be fully
faithful.  This problem already occurs for \(G=\Z\).  The issue is that the
tensor product of \(\Sch(G)\) with itself plays a crucial role.  If we work in
\(\Mod_\alg\bigl(\Sch(G)\bigr)\), we have to deal with
\(\Sch(G)\otimes\Sch(G)\), which appears quite intractable.  In
\(\Mod\bigl(\Sch(G)\bigr)\) we meet instead the much simpler completion
\(\Sch(G\times G)\) of this space.

Now it is time to explain briefly how we do analysis.  I am an advocate of
bornologies as opposed to topologies.  This means working with bounded subsets
and bounded maps instead of open subsets and continuous maps.  General
bornological vector spaces behave better than general topological vector
spaces for purposes of representation theory and homological algebra (see
\cites{Meyer:Smooth, Meyer:Embed_derived}).  The spaces that we shall use here
carry both a bornology and a topology, and both structures determine each
other.  Therefore, readers who are familiar with topological vector spaces may
be able to follow this article without learning much about bornologies.  We
explain some notions of bornological analysis along the way because they may
be unfamiliar to many readers.

We let \(\Mod(G)\) be the category of smooth representations of~\(G\) on
bornological vector spaces as in~\cite{Meyer:Smooth}.  The algebras
\(\Hecke(G)\) and \(\Sch(G)\) are bornological algebras in a natural way.  A
smooth representation \(\pi\colon G\to\Aut(V)\) on a bornological vector
space~\(V\) is called tempered if its integrated form extends to a bounded
algebra homomorphism \(\Sch(G)\to\Endo(V)\).  We may identify \(\Mod(G)\) with
the category of essential (or non-degenerate) bornological left modules over
\(\Hecke(G)\) (\cite{Meyer:Smooth}).  As our notation suggests, this
identifies the subcategory \(\Mod\bigl(\Sch(G)\bigr)\) with the category of
essential bornological modules over \(\Sch(G)\).  We turn \(\Mod(G)\) and
\(\Mod\bigl(\Sch(G)\bigr)\) into exact categories using the class of
extensions with a bounded linear section.  The exact category structure allows
us to form the derived categories \(\Der(G)\) and \(\Der\bigl(\Sch(G)\bigr)\)
as in \cite{Keller:Handbook}.  Actually, the passage to derived categories is
rather easy in both cases because our categories have enough projective and
injective objects.

Equipping a vector space with the finest possible bornology, we identify the
category of vector spaces with a full subcategory of the category of
bornological vector spaces.  Thus \(\Mod_\alg(G)\) becomes a full subcategory
of \(\Mod(G)\).  Moreover, this embedding maps projective objects again to
projective objects.  Therefore, the induced functor \(\Der_\alg(G)\to\Der(G)\)
is still fully faithful.  Our main theorem asserts that the canonical functor
\(\Der\bigl(\Sch(G)\bigr)\to\Der(G)\) is fully faithful as well.  The basic
technology for its proof is already contained in \cites{Meyer:Embed_derived,
  Meyer:Poly_comb}.

In \cite{Meyer:Embed_derived}, I define the category \(\Mod(A)\) of essential
modules and its derived category \(\Der(A)\) for a ``quasi-unital''
bornological algebra~\(A\) and extend some homological machinery to this
setting.  A morphism \(A\to B\) is called \emph{isocohomological} if the
induced functor \(\Der(B)\to\Der(A)\) is fully faithful.
\cite{Meyer:Embed_derived} gives several equivalent characterisations of
isocohomological morphisms.  The criterion that is most easy to verify is the
following: let \(P_\bullet\to A\) be a projective \(A\)\nbd{}bimodule
resolution of~\(A\); then \(A\to B\) is isocohomological if and only if
\(B\hot_A P_\bullet\hot_A B\) is a resolution of~\(B\) (in both cases,
resolution means that there is a bounded contracting homotopy).

The article~\cite{Meyer:Poly_comb} deals with the special case of the
embedding \(\C[G]\to\Sch_1(G)\) for a finitely generated discrete group~\(G\)
and a certain Schwartz algebra \(\Sch_1(G)\), which is defined by
\(\ell_1\)\nbd{}estimates.  The chain complex whose contractibility decides
whether this embedding is isocohomological turns out to be a \emph{coarse
  geometric} invariant of~\(G\).  That is, it depends only on the
quasi-isometry class of a word-length function on~\(G\).  If the group~\(G\)
admits a sufficiently nice combing, then I construct an explicit contracting
homotopy of this chain complex.  Thus \(\C[G]\to\Sch_1(G)\) is
isocohomological for such groups.

The argument for Schwartz algebras of reductive \(p\)-adic groups follows the
same pattern.  Let \(H\subseteq G\) be some compact open subgroup and let
\(X\defeq G/H\).  This is a discrete space which inherits a canonical coarse
geometric structure from~\(G\).  Since~\(G\) is reductive, it acts properly
and cocompactly on a Euclidean building, namely, its affine Bruhat-Tits
building.  Such buildings are CAT(0) spaces and hence combable.  Since~\(X\)
is coarsely equivalent to the building, it is combable as well.  Thus the
geometric condition of~\cite{Meyer:Poly_comb} is easily fulfilled for all
reductive \(p\)-adic groups.  However, we also have to check that the
constructions in~\cite{Meyer:Poly_comb} are compatible with uniform smoothness
of functions because~\(G\) is no longer discrete.  This forces us to look more
carefully at the geometry of the building.

\section{Bornological analysis}
\label{sec:bornologies}

Algebras like the Schwartz algebra \(\Sch(G)\) of a reductive \(p\)\nbd{}adic
group carry an additional structure that allows to do analysis in them.  The
homological algebra for modules over such algebras \emph{simplifies} if we
take this additional structure into account.  One reason is that the complete
tensor product \(\Sch(G)\hot\Sch(G)\) can be identified with \(\Sch(G^2)\)
(Lemma~\ref{lem:Sch_tensor}).

It is customary to describe this additional structure using a locally convex
topology.  We prefer to use bornologies instead.  This means that we work with
bounded subsets and bounded operators instead of open subsets and continuous
operators.  A basic reference on bornologies
is~\cite{Hogbe-Nlend:Bornologies}.  We use bornologies because of their
advantages in connection with homological algebra
(see~\cite{Meyer:Embed_derived}).

We mainly need bornological vector spaces that are complete and convex.
Therefore, we drop these adjectives and tacitly require all bornologies to be
complete and convex.  When we use incomplete bornologies, we explicitly say
so.

We need two classes of examples: fine bornologies and von Neumann bornologies.
Let~\(V\) be a vector space over~\(\C\).  The \emph{fine bornology}
\(\Fine(V)\) is the finest possible bornology on~\(V\).  A subset \(T\subseteq
V\) is bounded in \(\Fine(V)\) if and only if there is a finite-dimensional
subspace \(V_T\subseteq V\) such that~\(T\) is a bounded subset of
\(V_T\cong\R^n\) in the usual sense.  We also write \(\Fine(V)\) for~\(V\)
equipped with the fine bornology.

Any linear map \(\Fine(V)\to W\) is bounded.  This means that \(\Fine\) is a
fully faithful functor from the category of vector spaces to the category of
bornological vector spaces that is left-adjoint to the forgetful functor in
the opposite direction.

Let~\(V\) be a (quasi)complete locally convex topological vector space.  A
subset \(T\subseteq V\) is called \emph{von Neumann bounded} if it is absorbed
by all neighbourhoods of zero.  These subsets form a bornology on~\(V\) called
the \emph{von Neumann bornology} (following~\cite{Hogbe-Nlend:Bornologies}).
We write \(\Vonn(V)\) for~\(V\) equipped with this bornology.

This defines a functor \(\Vonn\) from topological to bornological vector
spaces.  Its restriction to the full subcategory of Fréchet spaces or, more
generally, of LF-spaces, is fully faithful.  That is, a linear map between
such spaces is bounded if and only if it is continuous.  A crucial advantage
of bornologies is that joint boundedness is much weaker than joint continuity
for multilinear maps: if \(V_1,\dotsc,V_n,W\) are (quasi)complete locally
convex topological vector spaces, then any separately continuous \(n\)-linear
map \(V_1\times\dotsm\times V_n\to W\) is (jointly) bounded.  The converse
also holds under mild hypotheses.

Let~\(G\) be a reductive \(p\)-adic group.  We carefully explain how the
Schwartz algebra \(\Sch(G)\) looks like as a bornological algebra.  The most
convenient definition for our purposes is due to Marie-France Vignéras
(\cite{Vigneras:Dimension}).  Let \(\sigma\colon G\to\N\) be the usual scale
on~\(G\).  It can be defined using a representation of~\(G\).  Let \(L_2(G)\)
be the Hilbert space of square-integrable functions with respect to some Haar
measure on~\(G\).  Let
\[
L_2^\sigma(G) \defeq
\{f\colon G\to\C\mid \text{\(f\cdot \sigma^k\in L_2(G)\) for all \(k\in\N\)}\}.
\]
A subset \(T\subseteq L_2^\sigma(G)\) is bounded if for all \(k\in\N\) there
exists a constant \(C_k\in\R_+\) such that \(\norm{f\cdot\sigma^k}_{L_2(G)}\le
C_k\) for all \(f\in T\).  This is the von Neumann bornology with respect to
the Fréchet topology on \(L_2^\sigma(G)\) defined by the sequence of
semi-norms
\[
\norm{f}^k_2\defeq \norm{f\cdot\sigma^k}_{L_2(G)}.
\]

Let \(\COM(G)\) be the set of compact open subgroups of~\(G\), ordered by
inclusion.  For \(U\in\COM(G)\), let \(\Sch(G\biinv U)=L_2^\sigma(G\biinv
U)\) be the subspace of \(U\)\nbd{}bi-invariant functions in \(L_2^\sigma(G)\).
We give \(\Sch(G\biinv U)\) the subspace bornology, that is, a subset is
bounded if and only if it is bounded in \(L_2^\sigma(G)\).  Finally, we let
\[
\Sch(G) \defeq \varinjlim \Sch(G\biinv U),
\]
where~\(U\) runs through \(\COM(G)\).  We equip \(\Sch(G)\) with the
direct-limit bornology.  That is, a subset of \(\Sch(G)\) is bounded if and
only if it is a bounded subset of \(\Sch(G\biinv U)\) for some
\(U\in\COM(G)\).  We may also characterise this bornology as the von Neumann
bornology with respect to the direct-limit topology on \(\Sch(G)\), using the
well-known description of bounded subsets in LF-spaces (see
\cite{Treves:Kernels}*{Proposition 14.6}).

\begin{lemma}[\cite{Vigneras:Dimension}]
  \label{lem:Schwartz_compare}
  The definition of the Schwartz algebra above agrees with the one of
  Harish-Chandra in \cites{Silberger:Reductive_padic,
    Waldspurger:Plancherel_formula}.
\end{lemma}

\begin{proof}
  The first crucial point is that the space of \emph{double} cosets \(G\biinv
  U\)\mdash as opposed to the group~\(G\) itself\mdash has polynomial growth
  with respect to the scale~\(\sigma\).  It suffices to check this for a good
  maximal compact subgroup~\(K\) because the map \(G\biinv U\to G\biinv K\)
  is finite-to-one.  By the Iwasawa decomposition, the double cosets in
  \(G\biinv K\) can be parametrised by points in a maximal split torus.  The
  scale on~\(G\) restricts to a standard word-length function on this torus,
  so that we get the desired polynomial growth.  As a result, there exists
  \(d>0\) such that \(\sum_{x\in G\biinv U} \sigma^{-d}(x)\) is bounded.

  Moreover, we need the following relationship between the growth of the
  double cosets \(UxU\) and the Harish-Chandra spherical function~\(\Xi\):
  there are constants \(C,r>0\) such that
  \[
  \vol(UxU) \le C\sigma(x)^r \cdot \Xi(x)^{-2},
  \qquad \Xi(x)^{-2} \le C\sigma(x)^r \cdot \vol(UxU).
  \]
  This follows from Equation I.1.(5) and Lemma II.1.1 in
  \cite{Waldspurger:Plancherel_formula}.  Hence
  \begin{multline*}
    \int_G \abs{f(x)}^2\sigma(x)^s\,dx
    = \sum_{x\in G\biinv U} \abs{f(UxU)}^2 \sigma(x)^s \vol(UxU)
    \\ \le \sum_{x\in G\biinv U} \abs{f(UxU)}^2 \Xi(x)^{-2} C\sigma(x)^{r+s}
    \\ \le \max_{x\in G} {} \abs{f(x)}^2 \Xi(x)^{-2} \sigma(x)^{r+s+d} \sum_{y\in
      G\biinv U} C\sigma^{-d}(y).
  \end{multline*}
  A similar computation shows
  \[
  \int_G \abs{f(x)}^2\sigma(x)^s\,dx
  \ge \max_{x\in G} {} \abs{f(x)}^2 \Xi(x)^{-2} C^{-1}\sigma(x)^{s-r}.
  \]
  Therefore, the sequences of semi-norms \(\norm{f\sigma^s}_2\) and
  \(\norm{f\Xi^{-1}\sigma^s}_\infty\) for \(s\in\N\) are equivalent and define
  the same function space \(\Sch(G\biinv U)\).
\end{proof}

Convolution defines a continuous bilinear map \(\Sch(G\biinv U) \times
\Sch(G\biinv U)\to\Sch(G\biinv U)\) for any \(U\in\COM(G)\) by
\cite{Waldspurger:Plancherel_formula}*{Lemme III.6.1}.  Since \(\Sch(G\biinv
U)\) is a Fréchet space, boundedness and continuity of the convolution are
equivalent.  Since any bounded subset of \(\Sch(G)\) is already contained in
\(\Sch(G\biinv U)\) for some~\(U\), the convolution is a bounded bilinear map
on \(\Sch(G)\), so that \(\Sch(G)\) is a bornological algebra.  In contrast,
the convolution on \(\Sch(G)\) is only separately continuous.

Now we return to the general theory and define the \(\Hom\) functor and the
tensor product.  Let \(\Hom(V,W)\) be the vector space of bounded linear maps
\(V\to W\).  A subset~\(T\) of \(\Hom(V,W)\) is bounded if and only if it is
\emph{equibounded}, that is, \(\{f(v)\mid f\in T,\ v\in S\}\) is bounded for
any bounded subset \(S\subseteq V\).  This bornology is automatically complete
if~\(W\) is.

The \emph{complete projective bornological tensor product}~\(\hot\) is defined
in~\cite{Hogbe-Nlend:Completions} by the expected universal property: it is a
bornological vector space \(V\hot W\) together with a bounded bilinear map
\(b\colon V\times W\to V\hot W\) such that \(l\mapsto l\circ b\) is a
bijection between bounded linear maps \(V\hot W\to X\) and bounded bilinear
maps \(V\times W\to X\).  This tensor product enjoys many useful properties.
It is commutative, associative, and commutes with direct limits.  It satisfies
the adjoint associativity relation
\begin{equation}
  \label{eq:adjoint_associativity}
  \Hom(V\hot W,X)\cong \Hom\bigl(V,\Hom(W,X)\bigr).
\end{equation}
Therefore, a bornological module over a bornological algebra~\(A\) can be
defined in three equivalent ways, using a bounded linear map \(A\to\Endo(V)\),
a bounded bilinear map \(A\times V\to V\), or a bounded linear map \(A\hot
V\to V\).

Let~\(\otimes\) be the usual tensor product of vector spaces.  The fine
bornology functor is compatible with tensor products; that is, the obvious map
\(V\otimes W\to V\hot W\) is a bornological isomorphism
\begin{equation}
  \label{eq:tensor_fine}
  \Fine(V\otimes W) \cong \Fine(V)\hot\Fine(W)
\end{equation}
for any two vector spaces \(V\) and~\(W\).  More generally, if~\(W\) is any
bornological vector space, then the underlying vector space of \(\Fine(V)\hot
W\) is equal to the purely algebraic tensor product \(V\otimes W\).  A subset
\(T\subseteq V\otimes W\) is bounded if and only if there is a
finite-dimensional subspace \(V_T\subseteq V\) such that~\(T\) is contained in
and bounded in \(V_T\otimes W\cong \R^n\otimes W\cong W^n\).  Here~\(W^n\)
carries the direct-sum bornology.  The reason for this is that~\(\hot\)
commutes with direct limits.

If \(V_1\) and~\(V_2\) are Fréchet-Montel spaces, then we have a natural
isomorphism
\begin{equation}
  \label{eq:tensor_Frechet}
  \Vonn(V_1\hot_\pi V_2) \cong \Vonn(V_1)\hot\Vonn(V_2),
\end{equation}
where~\(\hot_\pi\) denotes the \emph{complete projective topological tensor
  product} (see \cites{Grothendieck:Produits, Treves:Kernels}).  This
isomorphism is proven in \cite{Meyer:Analytic}*{Appendix A.1.4}, based on
results of Alexander Grothendieck.  The Montel condition means that all von
Neumann bounded subsets are precompact (equivalently, relatively compact).

\begin{lemma}
  \label{lem:Sch_tensor}
  Let~\(G\) be a reductive \(p\)-adic group.  Then
  \(\Sch(G)\hot\Sch(G)\cong\Sch(G\times G)\).
\end{lemma}

\begin{proof}
  It is shown in~\cite{Vigneras:Dimension} that \(\Sch(G\biinv U)\) is a
  nuclear Fréchet space for all \(U\in\COM(G)\); in fact, this follows easily
  from the proof of Lemma~\ref{lem:Schwartz_compare}.  Equip~\(G^2\) with the
  scale \(\sigma(a,b)\defeq\sigma(a)\sigma(b)\) for all \(a,b\in G\).  By
  definition, \(\Sch(G^2) \cong \varinjlim \Sch(G^2\biinv U^2)\).
  Since~\(\hot\) commutes with direct limits,
  \[
  \Sch(G)\hot\Sch(G)
  \cong \varinjlim \Sch(G\biinv U)\hot \Sch(G\biinv U)
  \]
  as well.  It remains to prove \(\Sch(G\biinv U)^{\hot 2}\cong
  \Sch(G^2\biinv U^2)\).  Since these Fréchet spaces are nuclear, they are
  Montel spaces.  Hence \eqref{eq:tensor_Frechet} allows us to
  replace~\(\hot\) by~\(\hot_\pi\).  Now we merely have to recall the
  definition of nuclearity (see \cites{Grothendieck:Produits,
    Treves:Kernels}).
  
  Let \(V\) and~\(W\) be Fréchet spaces.  The natural map \(V\otimes
  W\to\Hom(V',W)\) defines another topology on \(V\otimes W\), which may be
  weaker than the projective tensor product topology.  A Fréchet space is
  nuclear if and only if this topology coincides with the projective tensor
  product topology.  Equivalently, there is only one topology on \(V\otimes
  W\) for which the canonical maps \(V\times W\to V\otimes W\) and \(V\otimes
  W\to \Hom(V',W)\) are continuous.  It is clear that the subspace topology
  from \(\Sch(G^2\biinv U^2)\) on \(\Sch(G\biinv U)\otimes \Sch(G\biinv
  U)\) has these two properties.  Hence it agrees with the projective tensor
  product topology.  Now the assertion follows because \(\Sch(G\biinv
  U)^{\otimes 2}\) is dense in \(\Sch(G^2\biinv U^2)\).
\end{proof}

\section{Basic homological algebra over the Hecke algebra}
\label{sec:homological_Hecke}

Throughout this section, \(G\) denotes a totally disconnected, locally compact
group, \(H\) denotes a fixed compact open subgroup of~\(G\), and \(X\defeq
G/H\).

Let~\(V\) be a bornological vector space and let \(\pi\colon G\to\Aut(V)\) be
a representation of~\(G\) by bounded linear operators.  The
representation~\(\pi\) is called \emph{smooth} if for any bounded subset
\(T\subseteq V\) there exists an open subgroup \(U\subseteq G\) such that
\(\pi(g,v)=v\) for all \(g\in U\), \(v\in T\) (see~\cite{Meyer:Smooth}).  For
example, the left and right regular representations of~\(G\) on \(\Sch(G)\)
are smooth.  If~\(V\) carries the fine bornology, the definition above is
equivalent to the usual notion of a smooth representation on a vector space.

Let \(\Mod(G)\) be the category of smooth representations of~\(G\) on
bornological vector spaces; its morphisms are the \(G\)\nbd{}equivariant
bounded linear maps.  Let \(\Mod_\alg(G)\) be the category of smooth
representations of~\(G\) on \(\C\)\nbd{}vector spaces.  The fine bornology
functor identifies \(\Mod_\alg(G)\) with a full subcategory of \(\Mod(G)\).

Let \(\Hecke(G)\) be the Hecke algebra of~\(G\); its elements are the locally
constant, compactly supported functions on~\(G\).  The convolution is defined
by
\[
f_1*f_2(g) = \int_G f_1(x)f_2(x^{-1}g)\,dx
\]
for some left-invariant Haar measure~\(dx\); we normalise it so that
\(\vol(H)=1\).  We equip \(\Hecke(G)\) with the fine bornology, so that
\(\Hecke(G)\in\Mod_\alg(G)\subseteq\Mod(G)\).  More generally, given any
bornological vector space~\(V\), we let \(\Hecke(G,V)\defeq \Hecke(G)\hot V\).
The underlying vector space of \(\Hecke(G,V)\) is just \(\Hecke(G)\otimes V\)
because \(\Hecke(G)\) carries the fine bornology.  Hence \(\Hecke(G,V)\) is
the space of locally constant, compactly supported functions \(G\to V\).  The
\emph{left regular representation}~\(\lambda\) and the \emph{right regular
  representation}~\(\rho\) of~\(G\) on \(\Hecke(G,V)\) are defined by
\[
\lambda_g f(x)\defeq f(g^{-1}x),
\qquad
\rho_g f(x)\defeq f(xg)
\]
as usual.  They are both smooth.

Any continuous representation \(\pi\colon G\to\Aut(V)\) on a bornological
vector space~\(V\) can be integrated to a bounded algebra homomorphism
\(\Hecke(G)\to\Endo(V)\), which we again denote by~\(\pi\).  By adjoint
associativity, this corresponds to a map
\[
\pi_*\colon \Hecke(G,V) = \Hecke(G)\hot V\to V,
\qquad f\mapsto \int_G \pi\bigl(g,f(g)\bigr)\,dg.
\]
The map~\(\pi_*\) is \(G\)\nbd{}equivariant if~\(G\) acts on \(\Hecke(G,V)\)
by~\(\lambda\).  By \cite{Meyer:Smooth}*{Proposition 4.7}, the
representation~\(\pi\) is smooth if and only if~\(\pi_*\) is a bornological
quotient map, that is, any bounded subset of~\(V\) is of the form \(\pi_*(T)\)
for some bounded subset \(T\subseteq \Hecke(G,V)\).  Even more, if~\(\pi\) is
smooth, then~\(\pi_*\) has a bounded linear section.  Namely, we can use
\begin{equation}
  \label{eq:sigma_H}
  \sigma_H\colon V\to\Hecke(G,V),
  \qquad \sigma_H v(g) = \pi(g^{-1},v) 1_H(g),
\end{equation}
where~\(1_H\) denotes the characteristic function of~\(H\).  Thus the category
\(\Mod(G)\) becomes isomorphic to the category \(\Mod\bigl(\Hecke(G)\bigr)\)
of essential modules over \(\Hecke(G)\) (see \cite{Meyer:Smooth}*{Theorem
  4.8}).  The term ``essential'' is a synonym for ``non-degenerate'' that is
not as widely used for other purposes.

Let~\(\Extensions\) be the class of all extensions in \(\Mod(G)\) that have a
bounded linear section.  This turns \(\Mod(G)\) into an exact category in the
sense of Daniel Quillen.  Hence the usual machinery of homological algebra
applies to \(\Mod(G)\): we can form a derived category \(\Der(G)\) and derived
functors (see \cites{Keller:Handbook, Meyer:Embed_derived}).  The exact
category \(\Mod(G)\) has enough projective and injective objects, so that the
usual recipes for computing derived functors apply.  We shall use the
following standard projective resolution in \(\Mod(G)\), which already occurs
in~\cite{Meyer:Poly_comb}.

The homogeneous space \(X\defeq G/H\) is discrete because~\(H\) is open
in~\(G\).  Let
\begin{equation}
  \label{eq:Xn_def}
  X_n \defeq
  \{(x_0,\dotsc,x_n)\in X^{n+1}\mid x_0\neq x_1,\dotsc, x_{n-1}\neq x_n\}.
\end{equation}
We equip \(\C[X_n]\) with the fine bornology.  We let~\(G\) act diagonally
on~\(X_n\) and equip \(\C[X_n]\) with the induced representation
\[
g\cdot f(x_0,\dotsc,x_n) \defeq f(g^{-1}x_0,\dotsc,g^{-1}x_n).
\]
The stabilisers of points
\[
\stab(x_0,\dotsc,x_n) = \bigcap_{j=0}^n x_jHx_j^{-1}
\]
are compact open subgroups of~\(G\) for each \((x_0,\dotsc,x_n)\in X_n\).  Let
\(X_n'\subseteq X_n\) be a subset that contains exactly one representative
from each orbit.  We get
\begin{equation}
  \label{eq:Xn_explicit}
  X_n = \coprod_{\xi\in X_n'} G/\stab(\xi),
  \qquad
  \C[X_n] = \bigoplus_{\xi\in X_n'} \C[G/\stab(\xi)].
\end{equation}
If \(U\in\COM(G)\), then we have a natural isomorphism
\[
\Hom_G(\C[G/U],V) \congto \Fix(U,V),
\qquad f\mapsto f(1_U).
\]
Since~\(U\) is compact, this is an exact functor of~\(V\), so that \(\C[G/U]\)
is a projective object of \(\Mod(G)\).  Therefore, \(\C[X_n]\) is projective
by~\eqref{eq:Xn_explicit}.

In the following, we view~\(X_n\) as a subset of \(\C[X_n]\) in the usual way.
We let \((x_0,\dotsc,x_n)=0\) if \(x_j=x_{j+1}\) for some
\(j\in\{0,\dotsc,n-1\}\).  Thus \((x_0,\dotsc,x_n)\in\C[X_n]\) is defined for
all \((x_0,\dotsc,x_n)\in X^{n+1}\).

We define the boundary map \(\delta=\delta_n\colon \C[X_{n+1}]\to\C[X_n]\) for
\(n\in\N\) by
\[
\delta\bigl( (x_0,\dotsc,x_{n+1}) \bigr) \defeq
\sum_{j=0}^{n+1} (-1)^j\cdot (x_0,\dotsc,\widehat{x_j},\dotsc,x_{n+1}),
\]
where~\(\widehat{x_j}\) means that~\(x_j\) is omitted.  In terms of functions,
we can write
\begin{equation}
  \label{eq:def_delta}
  \delta\phi(x_0,\dotsc,x_n) =
  \sum_{j=0}^{n+1} (-1)^j \sum_{y\in X}
  \phi(x_0,\dotsc,x_{j-1},y,x_j,\dotsc,x_n).
\end{equation}
The operators~\(\delta_n\) are \(G\)\nbd{}equivariant for all \(n\in\N\).  We
define the augmentation map \(\alpha\colon \C[X_0]\to\C\) by \(\alpha(x)=1\)
for all \(x\in X_0=X\).  It is \(G\)\nbd{}equivariant with respect to the
trivial representation of~\(G\) on~\(\C\).  It is easy to see that
\(\delta^2=0\) and \(\alpha\circ\delta_0=0\).  Hence we get a chain complex
\[
C_\bullet(X)\defeq (\C[X_n],\delta_n)_{n\in\N}
\]
over~\(\C\).  We also form the reduced complex \(\tilde{C}_\bullet(X)\), which
has \(\C[X_n]\) in degree \(n\ge1\) and \(\ker(\alpha\colon \C[X]\to\C)\) in
degree~\(0\).  The complex \(\tilde{C}_\bullet(X)\) is exact.  Thus
\(C_\bullet(X)\to\C\) is a projective resolution of~\(\C\) in
\(\Mod_\alg(G)\).

Next we define \emph{bivariant co-invariant spaces}.  For \(V,W\in\Mod(G)\),
let \(V\hot_G W\) be the quotient of \(V\hot W\) by the closed linear span of
\(v\otimes w- gv\otimes gw\) for \(v\in V\), \(w\in W\), \(g\in G\).  Thus
\(V\hot_G W\) is again a complete bornological vector space.  By definition,
we have
\[
V\hot_G W \cong W\hot_G V \cong (V\hot W)\hot_G \C,
\]
where we equip \(V\hot W\) with the diagonal representation and~\(\C\) with
the trivial representation of~\(G\).  If~\(X\) is another bornological vector
space, then we may identify \(\Hom(V\hot_G W,X)\) with the space of bounded
bilinear maps \(f\colon V\times W\to X\) that satisfy \(f(gv,gw)=f(v,w)\) for
all \(g\in G\), \(v\in V\), \(w\in W\).  This universal property characterises
\(V\hot_G W\) uniquely.  It follows from the defining property of~\(\hot\).

There is an alternative description of \(V\hot_G W\) in terms of
\(\Hecke(G)\)\nbd{}modules.  Turn \(V\) into a right and~\(W\) into a left
bornological \(\Hecke(G)\)\nbd{}module by
\[
v*f\defeq \int_G f(g) \,g^{-1} v\,dg,
\qquad
f*w\defeq \int_G f(g) \,g w\,dg
\]
for all \(f\in\Hecke(G)\), \(v\in V\), \(w\in W\).  Let \(V\hot_{\Hecke(G)}
W\) be the quotient of \(V\hot W\) by the closed linear span of \(v*f\otimes
w-v\otimes f*w\) for \(v\in V\), \(f\in\Hecke(G)\), \(w\in W\).

\begin{lemma}
  \label{lem:coinvariant}
  \(V\hot_G W=V\hot_{\Hecke(G)} W\), that is, the elements \(gv\otimes
  gw-v\otimes w\) and \(v*f\otimes w-v\otimes f*w\) generate the same closed
  linear subspace.
\end{lemma}

\begin{proof}
  We have to show \(\Hom(V\hot_G W,X)=\Hom(V\hot_{\Hecke(G)} W,X)\) for all
  bornological vector spaces~\(X\).  By definition, \(\Hom(V\hot_G W,X)\) is
  the space of bounded bilinear maps \(l\colon V\times W\to X\) that satisfy
  \(l(g^{-1}v,w)=l(v,gw)\) for all \(v\in V\), \(w\in W\), \(g\in G\).  This
  implies \(l(v*f,w)=l(v,f*w)\) for all \(v\in V\), \(w\in W\),
  \(f\in\Hecke(G)\).  Conversely, suppose \(l(v*f,w)=l(v,f*w)\).  Then
  \begin{multline*}
    l\bigl(gv,g\cdot (f*w)\bigr)
    = l(gv,(\delta_g*f)*w)
    \\ = l(gv*(\delta_g*f),w)
    = l\bigl((g^{-1}\cdot gv)*f,w\bigr)
    = l(v,f*w)
  \end{multline*}
  for all \(v\in V\), \(w\in W\), \(g\in G\), \(f\in\Hecke(G)\).  This implies
  \(l(gv,gw)=l(v,w)\) for all \(v\in V\), \(w\in W\), \(g\in G\) because any
  \(w\in W\) is fixed by some \(U\in \COM(G)\) and therefore of the form
  \(\mu_U*w\), where~\(\mu_U\) is the normalised Haar measure of~\(U\).
\end{proof}

Since~\(\hot_G\) is functorial in both variables, we can apply it to chain
complexes.  Especially, we get a chain complex of bornological vector spaces
\(V\hot_G C_\bullet(X)\).  Since \(C_\bullet(X)\) is a projective resolution
of the trivial representation, we denote the chain homotopy type of \(V\hot_G
C_\bullet(X)\) by \(V\Lhot_G \C\).  The homology vector spaces of \(V\Lhot_G
\C\) may be denoted \(\Tor^G_n(V,\C)\) or \(\Tor^{\Hecke(G)}_n(V,\C)\).
However, this passage to homology forgets an important part of the structure,
namely, the bornology.  Therefore, it is better to work with \(V\Lhot_G\C\)
instead.

If \(V\) and~\(W\) are just vector spaces, we can identify \(V\hot_G W\) with
a purely algebraic construction.  Let \(V\otimes_G W\) be the quotient of
\(V\otimes W\) by the linear span of \(v\otimes w- gv\otimes gw\) for \(v\in
V\), \(w\in W\), \(g\in G\).  Then
\[
\Fine(V\otimes_G W) \cong \Fine(V)\hot_G \Fine(W).
\]
This follows from~\eqref{eq:tensor_fine} and the fact that any linear subspace
of a fine bornological vector space is closed.  Therefore, \(\Fine(V)\hot_G
C_\bullet(X) \cong \Fine\bigl(V\otimes_G C_\bullet(X)\bigr)\), and
\(\Tor^G_*(V,\C)\) is the homology of the chain complex \(V\otimes_G
C_\bullet(X)\).  In this case, passage to homology is harmless because
\(V\Lhot_G\C\) carries the fine bornology; this implies that it is
quasi-isomorphic to its homology viewed as a complex with vanishing boundary
map.

Our next goal is to describe \(V\hot_G C_\bullet(X)\) (Proposition
\ref{pro:VhotGC_explicit}).  This requires some geometric preparations.

\begin{definition}
  \label{def:controlled_subset}
  Given a finite subset \(F\subseteq X\), we define the relation~\(\sim_F\)
  on~\(X\) by
  \begin{equation}
    \label{eq:closeness}
    x\sim_F y \iff (x,y)\in \bigcup_{g\in G} g\cdot (\{H\}\times F)
    \iff x^{-1}y\in HFH.
  \end{equation}
  Here we view \(x^{-1}y \in G\biinv H\) and \(HFH\subseteq G\biinv H\).
  
  A subset \(S\subseteq X_n\) is \emph{controlled by~\(F\)} if \(x_i\sim_F
  x_j\) for all \((x_0,\dotsc,x_n)\in S\) and all \(i,j\in\{0,\dotsc,n\}\).
  We call \(S\subseteq X_n\) \emph{controlled} if it is controlled by some
  finite~\(F\).  Roughly speaking, this means that all entries of~\(S\) are
  uniformly close.
\end{definition}

A subset \(S\subseteq X_n\) is controlled if and only if~\(S\) is
\(G\)\nbd{}finite, that is, there is a finite subset \(F\subseteq X_n\) such
that \(S\subseteq G\cdot F\).  This alternative characterisation will be used
frequently.  Definition~\ref{def:controlled_subset} emphasises a crucial link
between the controlled support condition and geometric group theory.

A \emph{coarse (geometric) structure} on a locally compact space such as~\(X\)
is a family of relations on~\(X\) satisfying some natural axioms due to John
Roe (see also~\cite{Emerson-Meyer:Dualizing}).  The subrelations of the
relations~\(\sim_F\) above define a coarse geometric structure on~\(X\) in
this sense.  Since it is generated by \(G\)\nbd{}invariant relations, it
renders the action of~\(G\) on~\(X\) isometric.  This property already
characterises the coarse structure uniquely: whenever a locally compact group
acts properly and cocompactly on a locally compact space, there is a unique
coarse structure for which this action is isometric (see
\cite{Emerson-Meyer:Dualizing}*{Example 6}).  Moreover, with this coarse
structure, the space~\(X\) is coarsely equivalent to~\(G\).

By definition, the notion of a controlled subset of~\(X_n\) depends only on
the coarse geometric structure of~\(X\).  Thus the space of functions
on~\(X_n\) of controlled support only depends on the large scale geometry
of~\(X\).

Although our main examples, reductive groups, are unimodular, we want to treat
groups with non-trivial modular function as well.  Therefore, we have to
decorate several formulas with modular functions.  We define the modular
homomorphism \(\Delta_G\colon G\to\R_{>0}\) by \(\Delta_G(g) d(g^{-1})=dg\)
and \(d(gh)=\Delta_G(h)\,dg\) for all \(h\in G\).

\begin{definition}
  \label{def:controlled}
  Let \(\Controlled(X_n,V)^\Delta\) be the space of all maps \(\phi\colon
  X_n\to V\) that have controlled support and satisfy the covariance condition
  \begin{equation}
    \label{eq:covariance}
    \phi(g\xi)=\Delta_G(g)^{-1} \pi\bigl(g,\phi(\xi)\bigr)
  \end{equation}
  for all \(\xi\in X_n\), \(g\in G\).  A subset \(T\subseteq
  \Controlled(X_n,V)^\Delta\) is bounded if \(\{\phi(\xi)\mid \phi\in T\}\) is
  bounded in~\(V\) for all \(\xi\in X_n\) and the supports of all \(\phi\in
  T\) are controlled by the same finite subset \(F\subseteq X\).
\end{definition}

\begin{proposition}
  \label{pro:VhotGC_explicit}
  For any \(V\in\Mod(G)\), there is a natural bornological isomorphism
  \[
  V\hot_G \C[X_n] \cong \Controlled(X_n,V)^\Delta.
  \]
  The induced boundary map on \(V\hot_G \C[X_n]\) corresponds to the boundary
  map
  \[
  \delta=\delta_n\colon
  \Controlled(X_{n+1},V)^\Delta\to\Controlled(X_n,V)^\Delta
  \]
  defined by~\eqref{eq:def_delta}.
\end{proposition}

We denote the resulting chain complex
\((\Controlled(X_n,V)^\Delta,\delta_n)_{n\in\N}\) by
\(\Controlled(X_\bullet,V)^\Delta\).

\begin{proof}
  The bifunctor~\(\hot_G\) commutes with direct limits and in particular with
  direct sums.  Hence~\eqref{eq:Xn_explicit} yields
  \begin{equation}
    \label{eq:orbit_decompose}
    V\hot_G \C[X_n]
    \cong \bigoplus_{\xi\in X_n'} V\hot_G \C[G\xi]
    \cong \bigoplus_{\xi\in X_n'} V\hot_G \C[G/\stab(\xi)].
  \end{equation}
  Fix \(\xi\in X_n'\) and let \(\Map(G\cdot\xi,V)^\Delta\) be the space of all
  maps from \(G\cdot\xi\) to~\(V\) that satisfy the covariance condition
  \eqref{eq:covariance}.  We equip \(\Map(G\cdot\xi,V)^\Delta\) with the
  product bornology as in Definition~\ref{def:controlled}.  We claim that the
  map
  \[
  I\colon V\hot\C[G\xi]\to\Map(G\xi,V),
  \qquad
  v\otimes\phi \mapsto
  [\eta \mapsto \int_G \pi(h^{-1},v)\cdot \phi(h\eta) \,dh],
  \]
  yields a bornological isomorphism
  \[
  V\hot_G\C[G\xi]\cong\Map_G(G\xi,V)^\Delta.
  \]
  We check that~\(I\) descends to \(V\hot_G \C[G\xi]\) and maps into
  \(\Map(G\xi,V)^\Delta\):
  \begin{multline*}
    I(gv\otimes g\phi)(\eta)
    = \int_G \pi(h^{-1},gv)\cdot g\phi(h\eta)\,dh
    = \int_G \pi(h^{-1}g,v)\cdot\phi(g^{-1}h\eta)\,dh
    \\ = \int_G \pi(h^{-1},v)\cdot\phi(h\eta)\,dh
    = I(v\otimes\phi),
  \end{multline*}
  \begin{multline*}
    I(v\otimes\phi)(g\eta)
    = \int_G \pi(h^{-1},v)\cdot \phi(hg\eta) \,dh
    \\ = \int_G \pi(gh^{-1},v)\phi(h\eta) \,d(hg^{-1})
    = \Delta_G(g^{-1}) \pi\bigl(g, I(v\otimes\phi)(\eta)\bigr).
  \end{multline*}
  Thus we get a well-defined map \(V\hot_G\C[G\xi]\to\Map(G\xi,V)^\Delta\).
  Evaluation at~\(\xi\) defines a bornological isomorphism
  \(\Map(G\xi,V)^\Delta\cong\Fix(\stab(\xi),V)\).  We claim that the latter is
  isomorphic to \(V\hot_G \C[G/\stab(\xi)]\).  Since \(\stab(\xi)\) is compact
  and open, the Haar measure \(\mu_{\stab(\xi)}\) of \(\stab(\xi)\) is an
  element of \(\Hecke(G)\).  Convolution on the right with
  \(\mu_{\stab(\xi)}\) is an idempotent left module homomorphism on
  \(\Hecke(G)\), whose range is \(\C[G/\stab(\xi)]\).  Since \(V\hot_G
  \Hecke(G)\cong V\) for all~\(V\), additivity implies that \(V\hot_G
  \C[G/\stab(\xi)]\) is equal to the range of \(\mu_{\stab(\xi)}\) on~\(V\),
  that is, to \(\Fix(\stab(\xi),V)\).  Thus we obtain an isomorphism
  \(V\hot_G\C[G\xi]\cong\Map_G(G\xi,V)^\Delta\), which can easily be
  identified with the map~\(I\).
  
  Recall that a subset of~\(X_n\) is controlled if and only if it meets only
  finitely many \(G\)\nbd{}orbits.  Therefore, we get the counterpart
  \(\Controlled(X_n,V)^\Delta \cong \bigoplus_{\xi\in X_n'}
  \Map(G\xi,V)^\Delta\) to~\eqref{eq:orbit_decompose}.  We can piece our
  isomorphisms on orbits together to an isomorphism
  \[
  I\colon
  V\hot_G \C[X_n] \to \Controlled(X_n,V)^\Delta,
  \qquad
  v\otimes\phi \mapsto [\xi \mapsto \int_G \pi(g^{-1},v) \phi(g\xi) \,dg].
  \]
  A straightforward computation yields \(\delta\circ
  I(v\otimes\phi)=I(v\otimes\delta\phi)\) for all \(v\in V\), \(\phi\in
  \C[X_n]\) with~\(\delta\) as in~\eqref{eq:def_delta}.  Therefore, \(I\)
  intertwines \(\ID\hot_G\delta\) and~\(\delta\).
\end{proof}

Now let \(C_\bullet(X,V)\defeq C_\bullet(X)\hot V\), equipped with the
diagonal representation of~\(G\).  Since \(C_\bullet(X)\) carries the fine
bornology, the underlying vector space of \(\C[X_n]\hot V\) may be identified
with the space of functions \(X_n\to V\) with finite support.

\begin{lemma}
  \label{lem:concrete_resolution}
  The chain complex \(C_\bullet(X,V)\) is a projective resolution of~\(V\) in
  \(\Mod(G)\).
\end{lemma}

\begin{proof}
  The complex \(C_\bullet(X,V)\) is exact because~\(\hot\) is exact on
  extensions with a bounded linear section.  We have \(\Hom_G(\C[G/U]\hot
  V,W)\cong \Hom_U(V,W)\) for any \(U\in\COM(G)\) and any smooth
  representation~\(W\).  Since this is an exact functor of~\(W\),
  \(\C[G/U]\hot V\) is projective.  Equation~\eqref{eq:Xn_explicit} shows that
  \(\C[X_n]\hot V\) is a direct sum of such representations and therefore
  projective as well.
\end{proof}

We view \(\Hecke(G)\) as a bimodule over itself in the usual way, by
convolution on the left and right.  Since right convolution commutes with the
left regular representation, the complex \(C_\bullet\bigl(X,\Hecke(G)\bigr)\)
is a complex of \(\Hecke(G)\)-bimodules.  The same reasoning as in the proof
of Lemma~\ref{lem:concrete_resolution} shows that it is a projective
\(\Hecke(G)\)-bimodule resolution of \(\Hecke(G)\).

For \(V,W\in\Mod(G)\), we let \(\Hom_G(V,W)\) be the space of bounded
\(G\)\nbd{}equivariant linear maps \(V\to W\), equipped with the equibounded
bornology.  It agrees with the space \(\Hom_{\Hecke(G)}(V,W)\) of bounded
linear \(\Hecke(G)\)\nbd{}module homomorphisms.  We also apply the bifunctor
\(\Hom_G\) to chain complexes.  In particular, we can plug in the projective
resolution \(C_\bullet(X,V)\) of Lemma~\ref{lem:concrete_resolution}.  The
homotopy type of the resulting cochain complex of bornological vector spaces
\(\Hom_G(C_\bullet(X,V),W)\) is denoted by \(\Right\Hom_G(V,W)\).  Its \(n\)th
cohomology vector space is \(\Ext^n_G(V,W)\).  As with \(V\Lhot_G \C\), it is
preferable to retain the cochain complex itself.

If \(V\) and~\(W\) carry the fine bornology, then \(C_\bullet(X,V) =
C_\bullet(X)\otimes V\) with the fine bornology.  Therefore,
\(\Right\Hom_G(V,W)\) is equal to the space of all \(G\)\nbd{}equivariant
linear maps \(C_\bullet(X)\otimes V\to W\).  Hence the \(\Ext\) spaces above
agree with the purely algebraic \(\Ext\) spaces.  In more fancy language, the
embedding \(\Mod_\alg(G)\to\Mod(G)\) induces a fully faithful functor between
the derived categories \(\Der_\alg(G)\to\Der(G)\).  This allows us to apply
results proven using analysis in a purely algebraic context.

\section{Isocohomological smooth convolution algebras}
\label{sec:isocoh_convolution}

We introduce a class of convolution algebras on totally disconnected, locally
compact groups~\(G\).  These have the technical properties that allow us to
formulate the problem.  Then we examine the notion of an isocohomological
embedding and formulate a necessary and sufficient condition for
\(\Hecke(G)\to\T(G)\) to be isocohomological.  This criterion involves the
contractibility of a certain bornological chain complex, which is quite close
to the one that arises in~\cite{Meyer:Poly_comb}.

\subsection{Unconditional smooth convolution algebras with rapid decay}
\label{sec:convolution_algebras}

Let~\(G\) be a totally disconnected, locally compact group.  Let
\(\sigma\colon G\to \R_{\ge1}\) be a scale with the following properties:
\(\sigma(ab)\le\sigma(a)\sigma(b)\) and \(\sigma(a)=\sigma(a^{-1})\);
\(\sigma\) is \(U\)\nbd{}bi-invariant for some \(U\in\COM(G)\); the
map~\(\sigma\) is proper, that is, the subsets
\begin{equation}
  \label{eq:SR}
  B_R(G)\defeq \{g\in G \mid \sigma(g)\le R\}
\end{equation}
are compact for all \(R\ge1\).  The usual scale on a reductive \(p\)\nbd{}adic
groups has these properties.  If the group~\(G\) is finitely generated and
discrete, then \(\sigma=1+\ell\) or \(2^\ell\) for a word-length
function~\(\ell\) are good, inequivalent choices.

Let \(U\in\COM(G)\).  Given sets \(S,S'\) of functions \(G\biinv U\to\C\) we
say that \emph{\(S'\) dominates~\(S\)} if for any \(\phi\in S\) there exists
\(\phi'\in S'\) with \(\abs{\phi'(g)}\ge\abs{\phi(g)}\) for all \(g\in G\).

\begin{definition}
  \label{def:usrdc}
  Let \(\T(G)\) be a bornological vector space of functions \(\phi\colon
  G\to\C\).  We call \(\T(G)\) an \emph{unconditional smooth convolution
    algebra of rapid decay} if it satisfies the following conditions:
  \begin{enumerate}[\ref{def:usrdc}.1.]
  \item \(\T(G)\) contains \(\Hecke(G)\);

  \item \(\Hecke(G)\) is dense in \(\T(G)\);
    
  \item the convolution extends to a bounded bilinear map
    \(\T(G)\times\T(G)\to\T(G)\);
    
  \item \(\T(G)=\varinjlim \T(G\biinv U)\) as bornological vector spaces,
    where~\(U\) runs through \(\COM(G)\) and \(\T(G\biinv U)\) is the space
    of \(U\)\nbd{}bi-invariant functions in \(\T(G)\);
    
  \item if a set of functions \(G\biinv U\to\C\) is dominated by a bounded
    subset of \(\T(G\biinv U)\), then it is itself a bounded subset of
    \(\T(G\biinv U)\);

  \item \(M_\sigma\) is a bounded linear operator on \(\T(G)\).

  \end{enumerate}
  The first four conditions define a \emph{smooth convolution algebra}, the
  fifth condition means that the convolution algebra is \emph{unconditional},
  the last one means that it has \emph{rapid decay}.
\end{definition}

An example of such a convolution algebra is the Schwartz algebra of a
reductive \(p\)-adic group.

Let \(\T(G)\) be an unconditional smooth convolution algebra of rapid decay.
A representation \(\pi\colon G\to\Aut(V)\) is called \emph{\(\T(G)\)-tempered}
if its integrated form extends to a bounded algebra homomorphism
\(\T(G)\to\Endo(V)\) or, equivalently, to a bounded bilinear map \(\T(G)\times
V\to V\).  The density of \(\Hecke(G)\) in \(\T(G)\) implies that this
extension is unique once it exists.  Furthermore, \(G\)\nbd{}equivariant maps
are \(\T(G)\)-module homomorphisms.  Since the subalgebras \(\T(G\biinv U)\)
are unital, the algebra \(\T(G)\) is ``quasi-unital'' in the notation
of~\cite{Meyer:Embed_derived}, so that the category \(\Mod\bigl(\T(G)\bigr)\)
of essential bornological left \(\T(G)\)-modules is defined.  This category is
naturally isomorphic to the category of \(\T(G)\)-tempered smooth
representations of~\(G\) (see~\cite{Meyer:Embed_derived}).  Thus
\(\Mod\bigl(\T(G)\bigr)\) is a full subcategory of \(\Mod(G)\).

The following lemmas prove some technical properties of \(\T(G)\) that are
obvious in most examples, anyway.  Define \(P_R\colon \T(G)\to\Hecke(G)\) by
\(P_R\phi(x)=\phi(x)\) for \(x\in B_R(G)\) and \(P_R\phi(x)=0\) otherwise,
with \(B_R(G)\) as in~\eqref{eq:SR}.

\begin{lemma}
  \label{lem:approximation_property}
  \(\lim_{R\to\infty} P_R(\phi)=\phi\) uniformly for~\(\phi\) in a bounded
  subset of \(\T(G)\).
\end{lemma}

\begin{proof}
  If \(T\subseteq\T(G)\) is bounded, then \(T\subseteq\T(G\biinv U)\) for
  some \(U\in\COM(G)\).  Shrinking~\(U\) further, we achieve that the
  scale~\(\sigma\) is \(U\)\nbd{}bi-invariant.  We may further assume that
  \(\phi'\in T\) whenever \(\phi'\colon G\biinv U\to\C\) is dominated by some
  \(\phi\in T\) because \(\T(G)\) is unconditional.  Since~\(M_\sigma\) is
  bounded, the subset \(M_\sigma(T)\subseteq\T(G\biinv U)\) is bounded as
  well.  For any \(\phi\in T\), we have \(\abs{\phi-P_R\phi} \le R^{-1}
  \abs{M_\sigma\phi}\), so that \(\phi-P_R\phi\in R^{-1} M_\sigma(T)\).  This
  implies uniform convergence \(P_R(\phi)\to\phi\) for \(\phi\in T\).
\end{proof}

In the following, we briefly write
\[
\T(G^2)\defeq \T(G)\hot\T(G).
\]
Lemma~\ref{lem:Sch_tensor} justifies this notation for Schwartz algebras of
reductive groups.  In general, consider the bilinear maps
\[
\T(G)\times\T(G)\to\C,
\qquad (\phi_1,\phi_2)\mapsto \phi_1(x)\phi_2(y)
\]
for \((x,y)\in G^2\).  They extend to bounded linear functionals on
\(\T(G^2)\) and hence map \(\T(G^2)\) to a space of smooth functions
on~\(G^2\).

\begin{lemma}
  \label{lem:TG_tensor}
  This representation of \(\T(G^2)\) by functions on~\(G^2\) is faithful, that
  is, \(\phi\in\T(G^2)\) vanishes once \(\phi(x,y)=0\) for all \(x,y\in G\).
\end{lemma}

\begin{proof}
  The claim follows easily from Lemma~\ref{lem:approximation_property} (this
  is a well-known argument in connection with Grothendieck's Approximation
  Property).  If \(\phi(x,y)=0\) for all \(x,y\in G\), then also \((P_R\hot
  P_R) \phi(x,y)=0\) for all \(R\in\N\), \(x,y\in G\).  Since \(P_R\hot P_R
  (\phi)\in\Hecke(G^2)\), this implies \(P_R\hot P_R(\phi)=0\) for all
  \(R\in\N\).  Lemma~\ref{lem:approximation_property} implies that \(P_R\hot
  P_R\) converges towards the identity operator on \(\T(G)\hot\T(G)\).  This
  yields \(\phi=0\) as desired.
\end{proof}

Hence we may view \(\T(G^2)\) as a space of functions on~\(G^2\).  It is easy
to see that \(\T(G^2)\) is again a smooth convolution algebra on~\(G^2\).
Equip~\(G^2\) with the scale \(\sigma_2(a,b)\defeq \sigma(a)\sigma(b)\) for
\(a,b\in G\).  Then the operator \(M_{\sigma_2}= M_\sigma\hot M_\sigma\) is
bounded, that is, \(\T(G^2)\) also satisfies the rapid decay condition.
However, \(\T(G^2)\) need not be unconditional.  We \emph{assume \(\T(G^2)\)
  to be unconditional} in the following.  This is needed for the proof of our
main theorem.

Let \(\Tc(G)\) be \(\Hecke(G)\) equipped with the subspace bornology from
\(\T(G)\).  This bornology is incomplete, of course.  Similarly, we let
\(\Tc(G^2)\) be \(\Hecke(G^2)\) equipped with the subspace bornology from
\(\T(G^2)\).

\begin{lemma}
  \label{lem:Tc_completion}
  The completions of \(\Tc(G)\) and \(\Tc(G^2)\) are naturally isomorphic to
  \(\T(G)\) and \(\T(G^2)\).
\end{lemma}

\begin{proof}
  It suffices to prove this for \(\Tc(G)\).  We verify by hand that \(\T(G)\)
  satisfies the universal property that defines the completion of \(\Tc(G)\).
  Alternatively, we could use general characterisations of completions
  in~\cite{Meyer:Born_Top}*{Section 4}.  We must show that any bounded linear
  map \(f\colon \Tc(G)\to W\) into a complete bornological vector space~\(W\)
  extends uniquely to a bounded linear map on \(\T(G)\).  By
  Lemma~\ref{lem:approximation_property}, the sequence of operators
  \(P_R\colon \T(G)\to\Tc(G)\) converges uniformly on bounded subsets towards
  the identity map on \(\T(G)\).  Hence any bounded extension~\(\bar{f}\)
  of~\(f\) satisfies \(\bar{f}(\phi)= \lim_{R\to\infty} f\circ P_R(\phi)\) for
  all \(\phi\in\T(G)\).  Conversely, this prescription defines a bounded
  linear extension of~\(f\).
\end{proof}

\subsection{Isocohomological convolution algebras}
\label{sec:isocoh}

Let~\(A\) be a quasi-unital algebra such as \(\Hecke(G)\) or \(\T(G)\).
In~\cite{Meyer:Embed_derived} I define the exact category \(\Mod(A)\) of
essential bornological left \(A\)\nbd{}modules and its derived category
\(\Der(A)\).  A bounded algebra homomorphism \(f\colon A\to B\) between two
quasi-unital bornological algebras induces functors \(f^*\colon
\Mod(B)\to\Mod(A)\) and \(f^*\colon \Der(B)\to\Der(A)\).  Trivially, if~\(f\)
has dense range then \(f^*\colon \Mod(B)\to\Mod(A)\) is fully faithful.  We
call~\(f\) \emph{isocohomological} if \(f^*\colon \Der(B)\to\Der(A)\) is fully
faithful as well (\cite{Meyer:Embed_derived}).  We are interested in the
embedding \(\Hecke(G)\to\T(G)\).  If it is isocohomological, we briefly say
that \(\T(G)\) is \emph{isocohomological}.  The following conditions are
proven in \cite{Meyer:Embed_derived}*{Theorem 35} to be equivalent to
\(\T(G)\) being isocohomological:
\begin{itemize}
\item \(V\Lhot_{\T(G)} W \cong V\Lhot_G W\) for all
  \(V,W\in\Der\bigl(\T(G)\bigr)\) (recall that
  \(\Lhot_{\Hecke(G)}\cong\Lhot_G\));

\item \(\Right\Hom_{\T(G)}(V,W)\cong\Right\Hom_G(V,W)\) for all
  \(V,W\in\Der\bigl(\T(G)\bigr)\);

\item the functor \(f^*\colon \Der\bigl(\T(G)\bigr)\to\Der(G)\) is fully
  faithful;

\item \(\T(G)\Lhot_G V \cong V\) for all \(V\in\Mod\bigl(\T(G)\bigr)\);

\item \(\T(G)\Lhot_G\T(G)\cong \T(G)\).

\end{itemize}
The last condition tends to be the easiest one to verify in practice.  We will
formulate it more concretely below.  The signs ``\(\cong\)'' in these
statements mean isomorphism in the homotopy category of chain complexes of
bornological vector spaces.  This is stronger than an isomorphism of homology
groups.  As a consequence, we have \(\Tor^G_n(V,W)\cong\Tor^{\T(G)}_n(V,W)\)
and \(\Ext^n_G(V,W)\cong\Ext^n_{\T(G)}(V,W)\) for all
\(V,W\in\Mod\bigl(\T(G)\bigr)\) if \(\T(G)\) is isocohomological.

Notions equivalent to that of an isocohomological embedding have been defined
independently by several authors, as kindly pointed out to me by A.\ Yu.\
Pirkovskii (see \cite{Pirkovskii:Stably_flat} and the references given there).
We warn the reader that in categories of topological algebras some of the
conditions above are no longer equivalent.  Namely, the cohomological
conditions in terms of the derived category and \(\Right\Hom\) are weaker than
the homological conditions involving~\(\Lhot\).

We have seen in Section~\ref{sec:homological_Hecke} that \(V\hot_G W\cong
(V\hot W)\hot_G\C\) for all \(V,W\in\Mod(G)\), where \(V\hot W\) is equipped
with the diagonal representation of~\(G\).  Since \(V\hot W\) is projective if
\(V\) or~\(W\) is projective, this implies an isomorphism
\[
V\Lhot_G W \cong (V\hot W)\Lhot_G \C
\]
for all \(V,W\in\Der(G)\).  Thus \(\T(G)\) is isocohomological if and only if
\[
\bigl(\T(G)\hot\T(G)\bigr) \Lhot_G \C = \T(G^2) \Lhot_G \C \cong \T(G).
\]
Here we equip \(\T(G^2)\) with the diagonal representation of~\(G\), which is
given by
\[
g\cdot f(x,y) \defeq \Delta_G(g) f(xg,g^{-1}y)
\]
for all \(g\in G\), \(f\in\T(G^2)\), \(x,y\in G\) because the left and right
\(\Hecke(G)\)-module structures on \(\T(G)\) are the integrated forms of the
left regular representation~\(\lambda\) and the \emph{twisted} right regular
representation \(\rho\cdot\Delta_G\).  The convolution map
\[
\T(G^2) = \T(G)\hot\T(G) \overset{*}\to\T(G)
\]
descends to a bounded linear map \(\T(G^2)\hot_G \C\to\T(G)\).  The latter map
is a bornological isomorphism because \(A\hot_A A\cong A\) for any
quasi-unital bornological algebra by \cite{Meyer:Embed_derived}*{Proposition
  16}.  Moreover, the convolution map \(\T(G^2)\to\T(G)\) has a bounded linear
section, namely, the map \(\T(G)\to\Hecke\bigl(G,\T(G)\bigr)\subseteq\T(G^2)\)
defined in~\eqref{eq:sigma_H}.

We may use the projective resolution \(C_\bullet(X)\to\C\) to compute
\(\T(G^2)\Lhot_G \C\).  Proposition~\ref{pro:VhotGC_explicit} identifies
\(\T(G^2)\hot_G C_\bullet(X)\) with
\(\Controlled(X_\bullet,\T(G^2))^\Delta\).  We augment this chain
complex by the map
\begin{equation}
  \label{eq:augmentation_map}
  \alpha\colon
  \T(G^2)\hot_G \C[X_0] \xrightarrow{\ID\hot_G\alpha}
  \T(G^2)\hot_G \C \xrightarrow[\cong]{*} \T(G)
\end{equation}
We let \(\tilde\Controlled(X_\bullet,\T(G^2))^\Delta\) be the subcomplex of
\(\Controlled(X_\bullet,\T(G^2))^\Delta\) that we get if we replace
\(\Controlled(X_0,\T(G^2))^\Delta\) by
\(\tilde\Controlled(X_0,\T(G^2))^\Delta\defeq\ker\alpha\).

\begin{proposition}
  \label{pro:isocoh_convolution}
  \(\T(G)\) is isocohomological if and only if
  \(\tilde\Controlled(X_\bullet,\T(G^2))^\Delta\) has a bounded contracting
  homotopy.
\end{proposition}

\begin{proof}
  Our discussion of the convolution map implies that the augmentation map
  in~\eqref{eq:augmentation_map} is a surjection with a bounded linear
  section.  Hence it is a chain homotopy equivalence \(\T(G^2)\Lhot_G \C \to
  \T(G)\) if and only if its kernel \(\tilde\Controlled(X_0,\T(G^2))^\Delta\)
  is contractible.
\end{proof}

To give the reader an idea why the various characterisations of
isocohomological embeddings listed above are equivalent, we explain how the
contractibility of \(\tilde\Controlled(X_\bullet,\T(G^2))^\Delta\) yields
isomorphisms \(\Right\Hom_{\T(G)}(V,W)\cong\Right\Hom_G(V,W)\) for
\(V,W\in\Mod\bigl(\T(G)\bigr)\).  Almost the same argument yields
\(V\Lhot_{\T(G)} W\cong V\Lhot_G W\).  The extension to objects of the derived
categories is a mere formality.

The space \(\T(G^2)\) carries a \(\T(G)\)-bimodule structure via
\(f_1*(f_2\otimes f_3)*f_4 = (f_1*f_2)\otimes (f_3*f_4)\).  This structure
commutes with the inner conjugation action, so that \(P_\bullet\defeq
\T(G^2)\hot_G C_\bullet(X)\) becomes a chain complex of bornological
\(\T(G)\)-bimodules over \(\T(G)\).  As above, we can compute these spaces
explicitly:
\[
\T(G^2)\hot_G \C[X_n]
\cong \bigoplus_{\xi\in X_n'} \Fix\bigl(\stab \xi,\T(G^2)\bigr).
\]
It is not hard to see that \(\T(G^2)\) is a projective object of
\(\Mod\bigl(\T(G^2)\bigr)\).  That is, \(\T(G^2)\) is a projective bimodule.
Since the summands of \(\T(G^2)\hot_G \C[X_n]\) are all retracts of
\(\T(G^2)\), we conclude that \(\T(G^2)\hot_G \C[X_n]\) is a projective
\(\T(G)\)-bimodule.
  
Suppose now that~\(P_\bullet\) is a resolution of \(\T(G)\).  Then it is a
projective \(\T(G^2)\)-bimodule resolution.  Since \(\T(G)\) is projective as
a right module, the contracting homotopy of~\(P_\bullet\) can be improved to
consist of bounded right \(\T(G)\)-module homomorphisms.  Therefore,
\(P_\bullet\hot_{\T(G)} V\) is again a resolution of~\(V\).  Explicitly,
\[
P_n\hot_{\T(G)} V \cong \bigoplus_{\xi\in X_n'} \Fix(\stab(\xi),\T(G)\hot V)
\]
because \(\T(G)\hot_{\T(G)} V\cong V\).  The summands are retracts of the
projective left \(\T(G)\)-module \(\T(G)\hot V\).  Hence
\(P_\bullet\hot_{\T(G)} V\) is a projective left \(\T(G)\)-module resolution
of~\(V\).  We use it to compute
\[
\Right\Hom_{\T(G)}(V,W) = \Hom_{\T(G)}(P_\bullet\hot_{\T(G)} V,W).
\]
Let \(U\in\COM(G)\) act on \(\T(G)\hot V\) by
\(\Delta_G|_U\cdot\rho\otimes\pi=\rho\otimes\pi\).  If \(f\colon V\to W\) is
bounded and \(U\)\nbd{}equivariant, then \(\phi\otimes v\mapsto \phi*f(v)\)
defines a bounded \(G\)\nbd{}equivariant linear map \(\Fix(U,\T(G)\hot V)\to
W\).  One can show that this establishes a bornological isomorphism
\[
\Hom_{\T(G)}(\Fix(U,\T(G)\hot V),W) \cong \Hom_U(V,W).
\]
This yields a natural isomorphism
\[
\Hom_{\T(G)}(P_n\hot_{\T(G)} V,W) \cong
\bigoplus_{\xi\in X_n'} \Hom_{\stab(\xi)} (V,W).
\]
The right hand side no longer depends on \(\T(G)\)!  Thus
\(\Hom_G(C_\bullet(X,V),W)\) is isomorphic to the same complex, and
\(\Right\Hom_{\T(G)}(V,W)\cong\Right\Hom_G(V,W)\) as asserted.

Next we simplify the chain complex \(\Controlled(X_\bullet,\T(G^2))^\Delta\).
To \(\phi\in\Controlled(X_n,\T(G^2))^\Delta\) we associate a function
\(\phi_*\colon G\times G\times X_n\to\C\) by \(\phi_*(g,h,\xi)\defeq
\phi(\xi)(g,h)\).  This identifies \(\Controlled(X_n,\T(G^2))^\Delta\) with a
space of functions on \(G^2\times X_n\) by Lemma~\ref{lem:TG_tensor}.  More
precisely, we get the space of functions \(\phi\colon G^2\times X_n\to\C\)
with the following properties:
\begin{itemize}
\item \(\supp\phi\subseteq G^2\times S\) for some controlled subset
  \(S\subseteq X_n\);
  
\item the function \((a,b)\mapsto \phi(a,b,\xi)\) belongs to \(\T(G^2)\) for
  all \(\xi\in X_n\);
    
\item \(\phi(ag,g^{-1}b,g^{-1}\xi)=\phi(a,b,\xi)\) for all \(\xi\in X_n\),
  \(g,a,b\in G\) (the two modular functions cancel).

\end{itemize}
The last condition means that~\(\phi\) is determined by its restriction to
\(\{1\}\times G\times X_n\) by \(\phi(a,b,\xi)=\phi(1,ab,a\xi)\).  Thus we
identify \(\Controlled(X_n,\T(G^2))^\Delta\) with the following function space
on \(G\times X_n\):

\begin{definition}
  \label{def:fun_GX}
  Let \(\Controlled(G\times X_n,\T)\) be the space of all functions
  \(\phi\colon G\times X_n\to\C\) with the following properties:
  \begin{enumerate}[\ref{def:fun_GX}.1.]
  \item \(\supp\phi\subseteq G\times S\) for some controlled subset
    \(S\subseteq X_n\);
  
  \item the function \((a,b)\mapsto \phi(ab,a\xi)\) belongs to \(\T(G^2)\) for
    all \(\xi\in X_n\).

  \end{enumerate}
  A subset \(T\subseteq \Controlled(G\times X_n,\T)\) is bounded if there is a
  controlled subset \(S\subseteq G\) such that \(\supp\phi\subseteq G\times
  S\) for all \(\phi\in T\) and if for any \(\xi\in X_n\), the set of
  functions \((a,b)\mapsto\phi(ab,a\xi)\) for \(\phi\in T\) is bounded in
  \(\T(G^2)\).
\end{definition}

The boundary map~\(\delta\) on \(\Controlled(X_n,\T(G^2))^\Delta\) corresponds
to the boundary map
\begin{multline*}
  \delta \colon \Controlled(G\times X_{n+1},\T)\to\Controlled(G\times X_n,\T),
  \\
  \delta\phi(g,x_0,\dotsc,x_n) =
  \sum_{j=0}^{n+1} (-1)^j \sum_{y\in X}
  \phi(g,x_0,\dotsc,x_{j-1},y,x_j,\dotsc,x_n).
\end{multline*}
The augmentation map \(\Controlled(X,\T(G^2))^\Delta\to\T(G)\) corresponds to
\begin{equation}
  \label{eq:augmentation}
  \alpha\colon \Controlled(G\times X,\T) \to \T(G),
  \qquad
  \alpha\phi(g) = \sum_{x\in X} \phi(g,x).
\end{equation}
The proofs are easy computations, which we omit.  Let
\(\tilde\Controlled(G\times X_0,\T)\subseteq \Controlled(G\times X_0,\T)\) be
the kernel of~\(\alpha\) and let \(\tilde\Controlled(G\times X_\bullet,\T)\)
be the bornological chain complex that we get if we replace
\(\Controlled(G\times X_0,\T)\) by \(\tilde\Controlled(G\times X_0,\T)\).
Thus
\[
\tilde\Controlled(G\times X_\bullet,\T)
\cong \tilde\Controlled(X_\bullet,\T(G^2))^\Delta.
\]

Smooth functions of compact support automatically satisfy both conditions in
Definition~\ref{def:fun_GX}, so that \(\Hecke(G)\otimes \C[X_n] \subseteq
\Controlled(G\times X_n,\T)\).  These embeddings are compatible with the
boundary and augmentation maps.  Thus \(\Hecke(G)\otimes
\tilde{C}_\bullet(X)\) becomes a subcomplex of \(\tilde\Controlled(G\times
X_\bullet,\T)\).  We write \(\tilde\Controlled(G\times X_\bullet,\Tc)\) and
\(\Controlled(G\times X_\bullet,\Tc)\) for the chain complexes
\(\Hecke(G)\otimes \tilde{C}_\bullet(X)\) and \(\Hecke(G)\otimes
C_\bullet(X)\) equipped with the incomplete subspace bornologies from
\(\Controlled(G\times X_\bullet,\T)\).  The complex
\(\tilde\Controlled(G\times X_\bullet,\Tc)\) is contractible because
\(\tilde{C}_\bullet(X)\) is.  However, the obvious contracting homotopy is
unbounded.

\begin{lemma}
  \label{lem:isocoh_criterion}
  Suppose that there is a contracting homotopy~\(D\) for
  \(\tilde{C}_\bullet(X)\) such that \(\ID_{\Hecke(G)} \otimes D\) is bounded
  on \(\tilde\Controlled(G\times X_\bullet,\Tc)\).  Then \(\T(G)\) is
  isocohomological.
\end{lemma}

\begin{proof}
  We claim that \(\tilde\Controlled(G\times X_\bullet,\T)\) is the completion
  of \(\tilde\Controlled(G\times X_\bullet,\Tc)\).  Then
  \(\tilde\Controlled(G\times X_\bullet,\T) \cong
  \tilde\Controlled(X_\bullet,\T(G^2))^\Delta\) inherits a bounded contracting
  homotopy because completion is functorial.
  Proposition~\ref{pro:isocoh_convolution} yields that \(\T(G)\) is
  isocohomological.  It remains to prove the claim.  We do this by reducing
  the assertion to Lemma~\ref{lem:Tc_completion}.  Since
  \(\tilde\Controlled(G\times X_\bullet,\T)\) is a direct summand in
  \(\Controlled(G\times X_\bullet,\T)\), it suffices to prove that
  \(\Controlled(G\times X_\bullet,\T)\) is the completion of
  \(\Controlled(G\times X_\bullet,\Tc)\).  Recall that \(X'_\bullet\) denotes
  a subset of~\(X_\bullet\) containing one point from each \(G\)\nbd{}orbit.
  The decomposition of~\(X_\bullet\) into \(G\)\nbd{}orbits yields a
  direct-sum decomposition
  \begin{equation}
    \label{eq:orbit_decompose_GXn}
    \Controlled(G\times X_\bullet,\Tc)
    \cong \bigoplus_{\xi\in X_\bullet'} \Fix\bigl(\stab \xi,\Tc(G^2)\bigr),
  \end{equation}
  and a similar decomposition for \(\Controlled(G\times X_\bullet,\T)\).  Here
  direct sums are equipped with the canonical bornology: a subset is bounded
  if it is contained in and bounded in a finite sub-sum.  The reason
  for~\eqref{eq:orbit_decompose_GXn} is that a subset of~\(X_n\) is controlled
  if and only if it meets only finitely many \(G\)\nbd{}orbits.  Since
  completion commutes with direct sums, the assertion now follows from
  Lemma~\ref{lem:Tc_completion}.
\end{proof}

\section{Contracting homotopies constructed from combings}
\label{sec:contract_combing}

In order to apply Lemma~\ref{lem:isocoh_criterion}, we have to construct
contracting homotopies of \(\tilde{C}_\bullet(X)\).  For this we use the
geometric recipes of~\cite{Meyer:Poly_comb}.  The only ingredient is a
sequence of maps \(p_k\colon X\to X\) with certain properties.  We first
explain how such a sequence of maps gives rise to a contracting homotopy~\(D\)
of \(\tilde{C}_\bullet(X)\).  Then we formulate conditions on \((p_k)\) and
prove that they imply boundedness of~\(D\).

\subsection{A recipe for contracting homotopies}
\label{sec:contract_recipe}

The construction of \(C_\bullet(X)\) and \(\tilde{C}_\bullet(X)\) is natural:
a map \(f\colon X\to X\) induces a chain map \(f_*\colon C_\bullet(X)\to
C_\bullet(X)\) by \(f_*\bigl((x_0,\dotsc,x_n)\bigr) \defeq
(f(x_0),\dotsc,f(x_n)\bigr)\) or, equivalently,
\begin{equation}
  \label{eq:Cbullet_functorial}
  f_*\phi(x_0,\dotsc,x_n) =
  \sum_{y_j\colon f(y_j)=x_j} \phi(y_0,\dotsc,y_n).
\end{equation}

Since \(\alpha\circ f_*=\alpha\), this restricts to a chain map on
\(\tilde{C}_\bullet(X)\).  We have \(\ID_*=\ID\) and \((fg)_*=f_*g_*\).
Let~\(p_0\) be be the constant map \(x\mapsto H\) for all \(x\in X\).  We
claim that \((p_0)_*=0\) on \(\tilde{C}_\bullet(X)\).  On \(\C[X_n]\) for
\(n\ge1\) this is due to our convention that \((x_0,\dotsc,x_n)=0\) if
\(x_i=x_{i+1}\) for some~\(i\).  For \(\phi\in\C[X_0]\), we get \((p_0)_*\phi
= \alpha(\phi)\cdot (H)\), where \((H)\in\C[X]\) is the characteristic
function of \(H\in X\).  This implies the claim.

Given maps \(f,f'\colon X\to X\), we define operators \(D_j(f,f')\colon
\C[X_n]\to\C[X_{n+1}]\) for \(j\in\{0,\dotsc,n\}\) by
\begin{equation}
  \label{eq:Dj}
  D_j(f,f')\bigl((x_0,\dotsc,x_n)\bigr) \defeq
  \bigl(f(x_0),\dotsc,f(x_j),f'(x_j),\dotsc,f'(x_n)\bigr)
\end{equation}
and let \(D(f,f')\defeq \sum_{j=0}^n (-1)^j D_j(f,f')\).  It is checked
in~\cite{Meyer:Poly_comb} that
\[
[\delta,D(f,f')] \defeq \delta \circ D(f,f') + D(f,f')\circ \delta = f'_*-f_*.
\]
Thus the chain maps~\(f_*\) on \(\tilde{C}_\bullet(X)\) for \(f\colon X\to X\)
are all chain homotopic.  In particular, \(D(\ID,p_0)\) is a contracting
homotopy of \(\tilde{C}_\bullet(X)\) because \((p_0)_*=0\).

However, this trivial contracting homotopy does not work for Lemma
\ref{lem:isocoh_criterion}.  Instead, we use a sequence of maps
\((p_k)_{k\in\N}\) with~\(p_0\) as above and \(\lim_{k\to\infty} p_k =\ID\),
that is, for each \(x\in X\) there is \(k_0\in\N\) such that \(p_k(x)=x\) for
all \(k\ge k_0\).  We let
\[
D_k\defeq D(p_k,p_{k+1}),
\qquad
D_{jk}\defeq D_j(p_k,p_{k+1}).
\]
Observe that \(D_{jk}\) vanishes on the basis vector \((x_0,\dotsc,x_n)\)
unless \(p_k(x_j)\neq p_{k+1}(x_j)\).  Therefore, all but finitely many
summands of
\[
D \defeq \sum_{k=0}^\infty D_k
\]
vanish on any given basis vector.  Thus~\(D\) is a well-defined operator on
\(C_\bullet(X)\).

The operator~\(D\) is a contracting homotopy of \(\tilde{C}_\bullet(X)\)
because
\begin{multline*}
  [D,\delta]
  = \sum_{k=0}^\infty [D(p_k,p_{k+1}),\delta]
  = \sum_{k=0}^\infty (p_{k+1})_* - (p_k)_*
  = \lim_{k\to\infty} (p_k)_* - (p_0)_*
  = \ID.
\end{multline*}
To verify this computation, plug in a basis vector and use that all but
finitely many terms vanish.  This is the operator we want to use in Lemma
\ref{lem:isocoh_criterion}.

\subsection{Sufficient conditions for boundedness}
\label{sec:conditions_boundedness}

Construct~\(D\) as above and let \(D'\defeq \ID_{\Hecke(G)}\otimes D\).  We
want this to be a bounded operator on \(\Controlled(G\times X_\bullet,\Tc)\).
For this, we impose three further conditions on \((p_k)\).  First, \((p_k)\)
should be a \emph{combing}.  This notion comes from geometric group theory and
is already used in~\cite{Meyer:Poly_comb}.  It allows us to control the
support of \(D'\phi\) for~\(\phi\) with controlled support.  Secondly, the
combing \((p_k)\) should be \emph{smooth}.  This allows us to control the
smoothness of the functions \((a,b)\mapsto D'\phi(ab,a\xi)\) on~\(G^2\) for
\(\xi\in X_n\).  Only the third condition involves the convolution algebra
\(\T(G)\).  It asks for a certain sequence of operators \(\Tc(G)\to\Tc(G^2)\)
to be equibounded.

The smoothness condition is vacuous for discrete groups.  The third condition
is almost vacuous for \(\ell_1\)-Schwartz algebras of discrete groups.  Hence
these two conditions are not needed in~\cite{Meyer:Poly_comb}.  In our
application to reductive groups, we construct the operators \((p_k)\) using
the retraction of the affine Bruhat-Tits building of the group along geodesic
paths.  This is a combing because Euclidean buildings are CAT(0) spaces.  Its
smoothness amounts to the existence of congruence subgroups.  The third
condition follows easily from Lemma~\ref{lem:Sch_tensor}.

We now formulate the above conditions on \((p_k)\) in detail and state the
main result.  We use the relation~\(\sim_F\) for a finite subset \(F\subseteq
X=G/H\) defined in~\eqref{eq:closeness} by \(x\sim_F y \iff x^{-1}y\in HFH\).

\begin{definition}
  \label{def:combing}
  A sequence of maps \((p_k)_{k\in\N}\) as above is called a \emph{combing}
  of~\(X\) if it has the following additional two properties:
  \begin{enumerate}[\ref{def:combing}.1.]
  \item there is a finite subset \(F\subseteq X\) such that \(p_k(x)\sim_F
    p_{k+1}(x)\) for all \(k\in\N\), \(x\in X\);
    
  \item for any finite subset \(F\subseteq X\) there is a finite subset
    \(\bar{F}\subseteq X\) such that \(p_k(x)\sim_{\bar{F}} p_k(y)\) for all
    \(k\in\N\) and \(x,y\in X\) with \(x\sim_F y\).

  \end{enumerate}
  We say that the combing has \emph{polynomial growth} (with respect to the
  scale~\(\sigma\)) if the least~\(k_0\) such that \(p_k(gH)=gH\) for all
  \(k\ge k_0\) grows at most polynomially in \(\sigma(g)\).  (This definition
  of growth differs slightly from the one in~\cite{Meyer:Poly_comb}.)
\end{definition}

We may view the sequence \(\bigl(p_k(x)\bigr)\) as a path from~\(H\) to~\(x\).
The conditions on a combing mean that these paths do not jump too far in each
step and that nearby elements have nearby paths \(\bigl(p_k(x)\bigr)\).

\begin{definition}
  \label{def:smooth_combing}
  A combing \((p_k)_{k\in G}\) of \(G/H\) is called \emph{smooth} if it has
  the following two properties:
  \begin{enumerate}[\ref{def:smooth_combing}.1.]
  \item all maps~\(p_k\) are \(\tilde{H}\)\nbd{}equivariant for some open
    subgroup \(\tilde{H}\subseteq H\);

  \item for any \(U\in\COM(G)\), there exists \(V\in\COM(G)\) such that
    \(aVb\subseteq UabU\) for all \(a,b\in G\) with \(p_k(abH)=aH\).
  \end{enumerate}
\end{definition}

\begin{definition}
  \label{def:compatible_combing}
  Let \((p_k)\) be a smooth combing of polynomial growth.  Define
  \[
  R_k\colon \Tc(G)\to \Tc(G^2),
  \qquad
  R_k\phi(a,b)\defeq
  \begin{cases}
    \phi(ab) & \text{if \(p_k(abH)=aH\)};
    \\ 0 & \text{otherwise}.
  \end{cases}
  \]
  We say that \((p_k)\) is \emph{compatible with \(\T(G)\)} if the sequence of
  operators~\((R_k)\) is equibounded.
\end{definition}

\begin{theorem}
  \label{the:combing_contracts}
  Let~\(G\) be a totally disconnected, locally compact group and let \(\T(G)\)
  be an unconditional smooth convolution algebra of rapid decay on~\(G\).
  Suppose also that the function space \(\T(G^2)\) on~\(G^2\) is
  unconditional.  If \(G/H\) for some compact open subgroup \(H\subseteq G\)
  admits a smooth combing of polynomial growth that is compatible with
  \(\T(G)\), then \(\T(G)\) is isocohomological.
\end{theorem}

\subsection{Proof of Theorem~\ref{the:combing_contracts}}
\label{sec:proof}

\begin{lemma}
  \label{lem:combing_controlled}
  Suppose that \((p_k)\) is a combing.  Then for any controlled subset
  \(S\subseteq X_n\) there is a controlled subset \(\bar{S}\subseteq X_{n+1}\)
  such that \(\supp\phi\subseteq S\) implies \(\supp D(\phi)\subseteq
  \bar{S}\) for all \(\phi\in\C[X_n]\).
\end{lemma}

\begin{proof}
  Since~\(S\) is controlled, there is a finite subset \(F\subseteq X\) such
  that \(x_i\sim_F x_j\) for all \(i,j\in\{0,\dotsc,n\}\), \((x_i)\in S\).
  Since~\((p_k)\) is a combing, we can find a finite subset \(\bar{F}\subseteq
  X\) such that \(p_k(x_i) \sim_{\bar{F}} p_{k+1}(x_i)\) and \(p_k(x_i)
  \sim_{\bar{F}} p_k(x_j)\) for all \(k\in\N\) and all \((x_j)\in S\).  Let
  \(F'\defeq H\bar{F}H\bar{F}H\subseteq G/H\).  If \(x\sim_{\bar{F}}
  y\sim_{\bar{F}} z\), then \(x\sim_{F'} z\).  Hence
  \(\bigl(p_k(x_0),\dotsc,p_k(x_j),p_{k+1}(x_j),\dotsc,p_{k+1}(x_n)\bigr)\) is
  controlled by~\(F'\) for all \((x_i)\in S\).  This means that all summands
  in \(D(x_0,\dotsc,x_n)\) are controlled by~\(F'\).
\end{proof}

Hence \(D'\defeq \ID_{\Hecke(G)}\otimes D\) preserves controlled supports in
\(\Controlled(G\times X_\bullet,\Tc)\).  We have seen
in~\eqref{eq:orbit_decompose_GXn} that
\[
\Controlled(G\times X_\bullet,\Tc)
\cong \bigoplus_{\xi\in X_\bullet'} \Fix\bigl(\stab \xi,\Tc(G^2)\bigr).
\]
The isomorphism sends \(\phi\colon G\times X_n\to\C\) to the family of
functions \((\phi_\xi)_{\xi\in X_n'}\) defined by \(\phi_\xi(a,b) \defeq
\phi(ab,a\xi)\).  Thus we may describe any operator on \(\Controlled(G\times
X_\bullet,\Tc)\) by a block matrix.  In particular, we get
\(D'=(D'_{\xi\eta})_{\xi,\eta\in X_\bullet'}\) with certain operators
\[
D'_{\xi\eta}\colon \Fix\bigl(\stab \eta,\Tc(G^2)\bigr)
\to \Fix\bigl(\stab \xi,\Tc(G^2)\bigr).
\]
The fact that~\(D'\) preserves controlled supports means that for
fixed~\(\eta\) we have \(D'_{\xi\eta}=0\) for all but finitely many~\(\xi\).
Thus the whole operator~\(D'\) is bounded if and only if all its matrix
entries \(D'_{\xi\eta}\) are bounded.

For \(j,n\in\N\), \(n\ge j\), define \(p_{jk} \colon X_n\to X^{n+2}\) by
\[
p_{jk}\bigl((x_0,\dotsc,x_n)\bigr)
\defeq \bigl(p_k(x_0),\dotsc,p_k(x_j),p_{k+1}(x_j),\dotsc,p_{k+1}(x_n)\bigr).
\]
Then the operator \(D_{jk}\colon \C[X_n]\to\C[X_{n+1}]\) is given by
\[
D_{jk}\phi(\xi) = \sum_{\textstyle\eta\in p_{jk}^{-1}(\xi)} \phi(\eta).
\]
Let \(D'_{jk}=\ID_{\Hecke(G)}\otimes D_{jk}\) and let \(D'_{jk,\xi\eta}\) be
the matrix entries of \(D'_{jk}\) with respect to the
decomposition~\eqref{eq:orbit_decompose_GXn}.  Thus \(D'_{\xi\eta} =
\sum_{k\in\N} \sum_{j=0}^n (-1)^j D'_{jk,\xi\eta}\).  Writing
\(\phi_\xi(a,b)=\phi(ab,a\xi)\) and \(\phi_\eta(ag,g^{-1}b)=\phi(ab,ag\eta)\),
we get
\begin{multline}
  \label{eq:Dprime_jk}
  D'_{jk,\xi\eta}\psi(a,b)
  = \sum_{\textstyle\{g\in G/\stab(\eta) \mid p_{jk}(ag\eta)=a\xi\}} \psi(ag,g^{-1}b)
  \\ = \vol(\stab \eta)^{-1}
  \int_{\textstyle\{g\in G \mid p_{jk}(ag\eta)=a\xi\}} \psi(ag,g^{-1}b) \,dg
\end{multline}
for \(\psi\in\Fix\bigl(\stab \eta,\Tc(G^2)\bigr)\).  The right hand side
of~\eqref{eq:Dprime_jk} makes sense for arbitrary \(\psi\in\Tc(G^2)\) and
extends \(D'_{jk,\xi\eta}\) to an operator on \(\Tc(G^2)\).  Now we fix
\(\xi,\eta\) until further notice and sometimes omit them from our notation.

Let \(U\subseteq\stab(\eta)\) be an open subgroup and \(\phi\in\Tc(G\biinv U)\).
Let \(\mu_U\in\Hecke(G)\) be the normalised Haar measure of~\(U\), that is,
\(\supp\mu_U=U\) and \(\mu_U(g)=\vol(U)^{-1}\) for \(g\in U\).
Equation~\eqref{eq:Dprime_jk} yields
\begin{multline*}
  \vol(\stab\eta) D'_{jk,\xi\eta}(\phi\otimes\mu_U)(a,b)
  \\ = \vol(U)^{-1} \int_{\textstyle\{g\in bU \mid p_{jk}(ag\eta)=a\xi\}} \phi(ag) \,dg
  =
  \begin{cases}
    \phi(ab) & \text{if \(p_{jk}(ab\eta)=a\xi\);}
    \\
    0 & \text{otherwise.}
  \end{cases}
\end{multline*}
Let \(\chi_{j,\xi\eta}(a,b)\) be the number of \(k\in\N\) with
\(p_{jk}(ab\eta)=a\xi\) and let
\[
\chi(a,b) = \chi_{\xi\eta}(a,b)
\defeq \sum_{j=0}^n (-1)^j \chi_{j,\xi\eta}(a,b).
\]
These numbers are finite for any \(a,b\in G\) for the same reason that
guarantees that the sum defining~\(D\) is finite on each basis vector.  Define
\begin{displaymath}
  A=A_{\xi\eta}\colon \Tc(G)\to\Tc(G^2),
  \qquad A\phi(a,b) = \chi(a,b)\cdot \phi(ab).
\end{displaymath}
Our computation shows that \(A\phi = D'_{\xi\eta}(\phi\otimes\mu_U)\cdot
\vol(\stab \eta)\) if~\(U\) is an open subgroup of \(\stab \eta\) and
\(\phi\in\Tc(G\biinv U)\).

\begin{lemma}
  \label{lem:Dprime_jk_one}
  The operator \(D'_{\xi\eta}\) is bounded if and only if~\(A_{\xi\eta}\) is
  bounded.
\end{lemma}

\begin{proof}
  The boundedness of \(D'_{\xi\eta}\) implies that~\(A\) is bounded on
  \(\Tc(G\biinv U)\) for sufficiently small~\(U\) and hence on all of \(\Tc(G)\).
  Suppose conversely that~\(A\) is bounded.  We turn \(\Tc(G^2)\) into an
  (incomplete) bornological right \(\Tc(G)\)-module by
  \[
  (f_1\otimes f_2)*f_3\defeq f_1\otimes (f_2*f_3).
  \]
  This bilinear map \(\Tc(G^2)\times\Tc(G)\to\Tc(G^2)\) is bounded because the
  convolution in \(\T(G)\) is bounded.  The operators \(D'_{jk,\xi\eta}\) and
  hence \(D'_{\xi\eta}\) are \(\Tc(G)\)-module homomorphisms
  by~\eqref{eq:Dprime_jk}.  Let \(U\subseteq\stab(\eta)\) be open and
  \(\phi_1,\phi_2\in\Tc(G\biinv U)\).  Then
  \[
  D'_{\xi\eta} (\phi_1\otimes\phi_2)
  = D'_{\xi\eta} (\phi_1\otimes\mu_U)*\phi_2
  = \vol(\stab\eta)^{-1} A(\phi_1)*\phi_2.
  \]
  This implies the boundedness of \(D'_{\xi\eta}\) because~\(U\) is
  arbitrarily small and~\(A\) and the convolution \(\Tc(G^2)\to\Tc(G)\) are
  bounded.
\end{proof}

\begin{lemma}
  \label{lem:smooth_combing}
  If the combing is smooth, then for any \(U\in\COM(G)\) there is
  \(V\in\COM(G)\) such that~\(A\) maps \(\Tc(G\biinv U)\) into
  \(\Tc(G^2\biinv V^2)\).
\end{lemma}

\begin{proof}
  Let \(a,b\in G\).  Clearly, \(\chi(a,by)=\chi(a,b)\) for
  \(y\in\stab(\eta)\).  Since~\(p_k\) is \(\tilde{H}\)\nbd{}equivariant, so is
  \(p_{jk}\).  Hence \(\chi(ha,b)=\chi(a,b)\) for \(h\in\tilde{H}\).
  Therefore, \(A\phi(ua,b)=A\phi(a,b)=A\phi(a,bu)\) for \(\phi\in\Tc(G\biinv
  U)\) provided \(U\subseteq \tilde{H}\cap\stab(\eta)\).  Moreover, we have
  \(A\phi(a,xb) = A\phi(ax,b)\) if \(x\in\stab(\xi)\).
  
  We may assume that the zeroth components \(\eta_0\) and~\(\xi_0\) are~\(H\):
  any \(G\)\nbd{}orbit on~\(X_n\) has such a representative.  Then
  \(p_{jk}(ab\eta)=a\xi\) implies \(p_k(abH)=aH\).  By the definition of a
  smooth combing, there is \(V\in\COM(G)\) such that \(aVb\subseteq UabU\)
  whenever \(p_{jk}(ab\eta)=a\xi\) for some \(k\in\N\).  Hence
  \(\phi(avb)=\phi(ab)\) for \((a,b)\in\supp\chi\) and \(v\in V\).  We may
  shrink~\(V\) such that \(V\subseteq\stab(\xi)\cap U\).  Then
  \(\chi(a,vb)=\chi(a,b)=\chi(av,b)\) as well, so that~\(A\) maps
  \(\Tc(G\biinv U)\) to \(\Tc(G^2\biinv V^2)\).
\end{proof}

Definition~\ref{def:compatible_combing} requires the following sequence of
operators to be equibounded:
\[
R_k\colon\Tc(G)\to\Tc(G^2),
\qquad
R_k\phi(a,b)\defeq
\begin{cases}
  \phi(ab) & \text{if \(p_k(abH)=aH\);}
  \\ 0 & \text{otherwise.}
\end{cases}
\]
Hence \(R\defeq \sum_{k\in\N} (k+1)^{-2} R_k\) is bounded.  We have
\(R\phi(a,b)= \phi(ab) \chi'(a,b)\), where
\[
\chi'(a,b)=\sum_{\textstyle\{k\in\N\mid p_k(abH)=aH\}} (k+1)^{-2}.
\]
Now let \(S\subseteq\Tc(G)\) be bounded.  Then \(S\subseteq\Tc(G\biinv U)\)
for some \(U\in\COM(G)\).  We have already found \(V\in\COM(G)\) such that
\(A\phi\in\Tc(G^2\biinv V^2)\) for all \(\phi\in\Tc(G\biinv U)\).  Since
pointwise multiplication by the scale~\(\sigma\) is bounded on \(\Tc(G\biinv
V)\), it follows that the set of functions \(\sigma(a)^N\sigma(b)^N R(S)\) is
bounded in \(\Tc(G^2\biinv V^2)\).  We claim that \(\sigma(a)^N\sigma(b)^N
R(S)\) dominates \(A(S)\) for sufficiently large~\(N\).  Since \(\T(G^2)\) is
unconditional, this implies the boundedness of~\(A\).

Since~\(A\) and~\(R\) are multiplication operators, the claim follows if
\(\chi'\sigma^N\) dominates \(\chi_{j,\xi\eta}\) for any fixed \(j,\xi,\eta\).
This is what we are going to prove.  Let \(k_0(ab)\) be the least~\(k_0\) such
that \(p_k(ab\eta_j)=ab\eta_j\) for all \(k\ge k_0\).  The polynomial growth
of the combing implies that \(k_0(ab)\) is dominated by \(C\sigma(ab)^N\le
C\sigma(a)^N\sigma(b)^N\) for sufficiently large \(C,N\).  Since \(\xi_j\neq
\xi_{j+1}\), we have \(p_{jk}(ab\eta)\notin G\xi\) for \(k\ge k_0\).  Hence,
if \(p_{jk}(ab\eta)=a\xi\), then \(k<k_0\).  We choose the set of
representatives~\(X_\bullet'\) such that \(\xi_0=H\) for all \(\xi\in
X_\bullet'\).  Then \(p_{jk}(ab\eta)=a\xi\) implies \(p_k(abH)=aH\).
Therefore, for each summand~\(1\) in \(\chi_{j,\xi\eta}(a,b)\) there is a
summand \(1/(k+1)^2\ge 1/k_0^2\) in \(\chi'(a,b)\).  This yields the desired
estimate
\(
\chi_{j,\xi\eta}(a,b)\le k_0(a,b)^2 \chi'(a,b)
\le C^2\sigma(a)^{2N}\sigma(b)^{2N} \chi'(a,b)
\).
Thus the operators \(A_{\xi\eta}\) are bounded for all \(\xi,\eta\).  This
implies boundedness of \(D'_{\xi\eta}\) by Lemma~\ref{lem:Dprime_jk_one}.  By
Lemma~\ref{lem:combing_controlled}, it follows that~\(D'\) is bounded.
Finally, Lemma~\ref{lem:isocoh_criterion} yields that \(\T(G)\) is
isocohomological.  This finishes the proof of Theorem
\ref{the:combing_contracts}.

\section{A smooth combing for reductive \(p\)-adic groups}
\label{sec:combing_reductive}

The following theorem is the main goal of this section.  In addition, we prove
a variant (Theorem~\ref{the:Sch_chi_isocoh}) that deals with the subcategories
of \(\chi\)\nbd{}homogeneous representations for a character \(\chi\colon
C(G)\to\U(1)\) of the connected centre of~\(G\).

\begin{theorem}
  \label{the:Sch_isocohomological}
  The Schwartz algebra \(\Sch(G)\) of a reductive \(p\)-adic group~\(G\) is
  isocohomological.
\end{theorem}

\begin{proof}
  We are going to apply Theorem~\ref{the:combing_contracts}.  Let~\(\sigma\)
  be the standard scale as in the definition of the Schwartz algebra and let
  \(\T(G)\defeq\Sch(G)\).  The space \(\T(G^2)\) is defined as
  \(\T(G)\hot\T(G)\).  This notation is permitted because of
  Lemma~\ref{lem:Sch_tensor}, which identifies \(\Sch(G)\hot\Sch(G)\) with the
  Schwartz algebra of~\(G^2\) (which is again a reductive \(p\)-adic group).
  Clearly, the Schwartz algebras \(\Sch(G)\) and \(\Sch(G^2)\) are
  unconditional smooth convolution algebras of rapid decay.

  We let \(\BT=\BT(G)\) be the \emph{affine Bruhat-Tits building} of~\(G\), as
  defined in \cites{Tits:Corvallis, Bruhat-Tits:Groupes_reductifs}.  This is a
  Euclidean building on which~\(G\) acts isometrically, properly, and
  cocompactly.  Let \(C(G)\) be the connected centre of~\(G\), so that the
  quotient \(G/C(G)\) is semi-simple.  In Section~\ref{sec:applications}, we
  will also use the variant \(\BT\bigl(G/C(G)\bigr)\) of \(\BT(G)\), which we
  call the \emph{semi-simple affine Bruhat-Tits building} of~\(G\).

  Let~\(G_\circ\) be the connected component of~\(G\) as an algebraic group.
  Thus~\(G_\circ\) is a reductive group and~\(G\) is a finite extension
  of~\(G_\circ\).  Inspection of the definition in~\cite{Tits:Corvallis} shows
  that the buildings for \(G\) and~\(G_\circ\) are equal.  We remark that it
  is not hard to reduce the case of general reductive \(p\)\nbd{}adic groups
  to the special case of connected semi-simple groups, or even connected
  simple groups.  At first I followed this route myself.  Eventually, it
  turned out that this intermediate step is unnecessary because all arguments
  work directly in the generality we need.

  Let \(\xi_0\in\BT\), \(H\defeq\stab(\xi_0)\), and \(X\defeq G/H\).  We have
  \(H\in\COM(G)\), and~\(X\) may be identified with the discrete subset
  \(G\xi_0\subseteq\BT\).  We need a combing of~\(X\).  As a preparation, we
  construct a combing of~\(\BT\), using that Euclidean buildings are CAT(0)
  spaces, that is, have ``non-positive curvature'' (see
  \cites{Bridson-Haefliger, Bruhat-Tits:Groupes_reductifs}).  In particular,
  any two points in~\(\BT\) are joined by a unique geodesic.  For
  \(\xi\in\BT\), let
  \[
  p(\xi)\colon [0,d(\xi,\xi_0)]\to\BT,
  \qquad t\mapsto p_t(\xi),
  \]
  be the unit speed geodesic segment from \(\xi_0\) to~\(\xi\); extend this by
  \(p_t(\xi)\defeq\xi\) for \(t>d(\xi,\xi_0)\).  Restricting to \(t\in\N\), we
  get a sequence of maps \(p_k\colon \BT\to\BT\).

  \begin{lemma}
    \label{lem:combing_building}
    The maps~\(p_k\) for \(k\in\N\) form a combing of linear growth of
    \(\BT\).
  \end{lemma}
  
  This means \(d\bigl(p_k(\xi),p_k(\eta)\bigr)\le R\cdot d(\xi,\eta)+R\) and
  \(d\bigl(p_k(\xi),p_{k+1}(\xi)\bigr)\le R\) for all \(\xi,\eta\in\BT\),
  \(k\in\N\), for some \(R>0\).  Linear growth means that the least~\(k_0\)
  such that \(p_k(\xi)\) is constant for \(k\ge k_0\) grows at most linearly
  in \(d(\xi,\xi_0)\).

  \begin{proof}
    By construction, \(d\bigl(p_s(\xi),p_t(\xi)\bigr)\le \abs{t-s}\) for all
    \(s,t\in\R_+\), \(\xi\in\BT\), and \(p_s(\xi)\) is constant for \(s\ge
    d(\xi,\xi_0)\).  The lemma follows if we prove the following claim:
    \(d\bigl(p_t(\xi),p_t(\eta)\bigr)\le d(\xi,\eta)\) for all
    \(\xi,\eta\in\BT\), \(t\in\R_+\).
    
    Fix \(\xi,\eta\in\BT\) and \(t\in\R_+\).  We may assume \(d(\xi,\xi_0)\ge
    d(\eta,\xi_0)\) (otherwise exchange \(\xi\) and~\(\eta\)) and
    \(d(\xi,\xi_0)\ge t\) (otherwise \(p_t(\xi)=\xi\) and \(p_t(\eta)=\eta\)).
    Let~\(d^*\) be the usual flat Euclidean metric on~\(\R^2\).  Let \(\xi^*\)
    and~\(\eta^*\) be points in~\(\R^2\) with
    \[
    d^*(\xi^*,0)=d(\xi,\xi_0),
    \qquad
    d^*(\eta^*,0)=d(\eta,\xi_0),
    \qquad
    d^*(\xi^*,\eta^*)=d(\xi,\eta).
    \]
    The CAT(0) condition means that distances between points on the boundary
    of the geodesic triangle \((\xi,\eta,\xi_0)\) are dominated by the
    distances between the corresponding points in the comparison triangle
    \((\xi^*,\eta^*,0)\).  Here the point \(p_t(\xi)\) corresponds to the
    point \(p_t(\xi)^*\) on \([0,\xi^*]\) of distance~\(t\) from the origin.
    The point \(p_t(\eta)\) corresponds to the point \(p_t(\eta)^*\) of
    distance \(\min\{t,d(\eta,\xi_0)\}\) from the origin.
    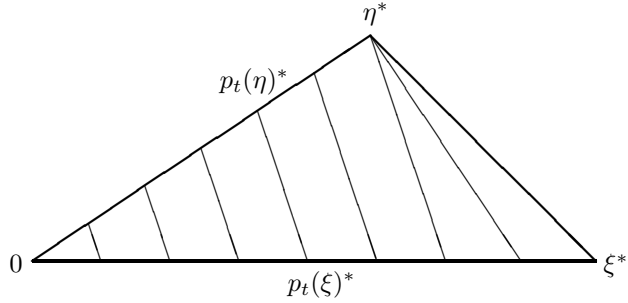
\begin{figure}[htb]
      \centering
      \setlength{\unitlength}{1cm}
      \begin{picture}(8,3.5)(0,0)
        \thicklines
        \put(0.2,0.2){\line(1,0){7.5}}
        \put(0.2,0.2){\line(3,2){4.5}}
        \put(4.7,3.2){\line(1,-1){3.0}}
        \put(-0.1,0.05){\(0\)}
        \put(7.8,0.05){\(\xi^*\)}
        \put(4.6,3.4){\(\eta^*\)}
        \thinlines
        \put(0.95,0.7){\line(1,-3){0.16667}}
        \put(1.7,1.2){\line(1,-3){0.33333}}
        \put(2.45,1.7){\line(1,-3){0.5}}
        \put(3.2,2.2){\line(1,-3){0.66667}}
        \put(3.95,2.7){\line(1,-3){0.83333}}
        \put(4.7,3.2){\line(1,-3){1}}
        \put(4.7,3.2){\line(2,-3){2}}
        \put(3.6,-0.2){\(p_t(\xi)^*\)}
        \put(2.7,2.5){\(p_t(\eta)^*\)}
      \end{picture}
      \caption{A comparison triangle}
      \label{fig:comparison_triangle}
    \end{figure}%
    An easy computation or a glance at Figure~\ref{fig:comparison_triangle}
    shows \(d^*(p_t(\xi)^*,p_t(\eta)^*) \le d^*(\xi^*,\eta^*)\).  By the
    CAT(0) condition, this implies \(d\bigl(p_t(\xi),p_t(\eta)\bigr) \le
    d(\xi,\eta)\).
  \end{proof}
  
  We identify \(G/H\) with the orbit \(G\xi_0\subseteq\BT\).  Since the group
  action is cocompact, there is some \(R>0\) such that for any \(\xi\in\BT\)
  there exists \(\xi'\in G\xi_0\) with \(d(\xi,\xi')<R\).  We let
  \(p_k'(\xi)\) for \(\xi\in\BT\) be a point in \(G\xi_0\) with
  \(d\bigl(p_k'(\xi),p_k(\xi)\bigr)< R\).  We claim that any such choice
  defines a combing of \(G/H\).
  
  If we equip \(G\xi_0\) with the metric~\(d\) from~\(\BT\), the maps~\(p_k'\)
  on \(G\xi_0\) still form a combing in the metric sense because they are
  ``close'' to the combing \((p_k)\).  The subspace metric from~\(\BT\) and
  the relations~\(\sim_F\) in Definition~\ref{def:combing} generate the same
  coarse geometric structure on~\(X\).  That is, for any \(R>0\) there is a
  finite subset \(F\subseteq X\) such that \(d(x,y)\le R\) implies \(x\sim_F
  y\), and for any finite subset \(F\subseteq X\) there is \(R>0\) such that
  \(x\sim_F y\) implies \(d(x,y)\le R\).  This is easy to verify by hand.
  Alternatively, it follows from the uniqueness of coarse structures mentioned
  after Definition~\ref{def:combing}.  Hence \((p_k')\) is a combing of
  \(G/H\) in the sense we need.
  
  We also need the combing to be smooth.  To get this, we must choose the base
  point~\(\xi_0\) and the approximations \(p_k'(\xi)\) more carefully.  This
  requires some geometric facts about the apartments in the building.
  Let~\(K\) be the non-Archimedean local field over which~\(G\) is defined.
  Let \(S\subseteq G\) be a maximal \(K\)\nbd{}split torus of~\(G\).  We do
  not distinguish in our notation between the algebraic groups \(S\) and~\(G\)
  and their locally compact groups of \(K\)\nbd{}rational points.  Let
  \(X^*(S)\) and \(X_*(S)\) be the groups of algebraic characters and
  cocharacters of~\(S\), respectively.  The \(\R\)\nbd{}vector space \(A\defeq
  X_*(S)\otimes\R\) is the \emph{basic apartment} of \((G,S)\).

  Let \(\Phi\subseteq X^*(S)\) be the set of roots of~\(G\) relative to~\(S\).
  Choose a simple system of roots \(\Delta\subseteq\Phi\) and let
  \(A_+\subseteq A\) be the corresponding closed Weyl chamber:
  \[
  A_+ \defeq
  \{x\in A\mid \text{\(\alpha(x)\ge0\) for all \(\alpha\in\Delta\)}\}.
  \]
  Let~\(W\) be the \emph{Weyl group} of the root system~\(\Phi\).  It is the
  Coxeter group generated by orthogonal reflections in the hyperplanes
  \(\alpha(x)=0\) for \(\alpha\in\Delta\).  The positive cone~\(A_+\) is a
  fundamental domain for this action, that is, \(W(A_+)=A\)
  (see~\cite{Humphreys:Coxeter_groups}).
  
  Let \(Z\subseteq G\) be the centraliser of~\(S\).  There is a canonical
  homomorphism \(\nu\colon Z\to A\) (see \cite{Tits:Corvallis}*{(1.2)}).  Its
  kernel is compact and its range is a lattice \(\Lambda\subseteq A\); that
  is, \(\Lambda\) is a discrete and cocompact subgroup of~\(A\).  Moreover,
  let \(\Lambda_+\defeq \Lambda\cap A_+\).  Since~\(\Lambda\) is free Abelian,
  we can lift it to a subgroup of~\(Z\) and view \(\Lambda\subseteq Z\subseteq
  G\).  Let \(\Phi_{\mathrm{af}}\) be the set of \emph{affine roots} as in
  \cite{Tits:Corvallis}*{(1.6)}.  These are affine functions \(\alpha\colon
  A\to\R\) of the form \(\alpha(x)=\alpha_0(x)+\gamma\) with
  \(\alpha_0\in\Phi\) and certain \(\gamma\in\R\).  Recall that
  \(\Phi\subseteq\Phi_{\mathrm{af}}\) and that \(\Phi_{\mathrm{af}}\) is
  invariant under translation by~\(\Lambda\).  The subsets of~\(A\) of the
  form \(\{x\in A\mid \alpha(x)=0\}\) and \(\{x\in A\mid \alpha(x)\ge0\}\) for
  \(\alpha\in\Phi_{\mathrm{af}}\) are called \emph{walls} and
  \emph{half-apartments}, respectively.

  We define the \emph{closure} \(\clos(\Omega)\subseteq A\) of a non-empty
  subset \(\Omega\subseteq A\) as the intersection of all closed
  half-apartments containing~\(\Omega\) (see
  \cite{Bruhat-Tits:Groupes_reductifs}*{(7.1.2)}).  We claim that
  \begin{equation}
    \label{eq:closure_twopoint}
    \clos(\{0,\xi\}) = A_+ \cap (\xi-A_+)
  \end{equation}
  for all \(\xi\in\Lambda_+\).  Let \(B\defeq A_+ \cap (\xi-A_+)\).  This is
  an intersection of half-apartments containing \(0\) and~\(\xi\) because
  \(\Phi\cup (\xi-\Phi)\subseteq\Phi_{\mathrm{af}}\).  Hence
  \(\clos(\{0,\xi\})\subseteq B\).  It remains to show that any
  half-apartment~\(C\) containing \(0\) and~\(\xi\) also contains~\(B\).
  Let~\(C\) be defined by the equation \(\alpha\ge0\) for some affine
  root~\(\alpha\) with linear part \(\alpha_0\in\Phi\).  We distinguish the
  cases \(\alpha_0>0\) and \(\alpha_0<0\).  If \(\alpha_0>0\), then
  \(\alpha(x)=\alpha_0(x)+\alpha(0)\ge\alpha_0(x)\) is non-negative
  on~\(A_+\); if \(\alpha_0<0\), then
  \(\alpha(\xi-x)=-\alpha_0(x)+\alpha(\xi)\ge-\alpha_0(x)\) is non-negative
  on~\(A_+\), so that~\(\alpha\) is non-negative on \(\xi-A_+\).  Thus
  \(B\subseteq C\) in either case.  This finishes the proof that
  \(\clos(\{0,\xi\})=B\).

  \begin{lemma}
    \label{lem:approx_good}
    There is \(R>0\) such that for all \(\xi\in\Lambda_+\) and all
    \(\eta\in\clos(\{0,\xi\})\), there is \(\eta'\in
    \Lambda\cap\clos(\{0,\xi\})\) with \(d(\eta,\eta')\le R\).
  \end{lemma}

  \begin{proof}
    Let \(\Delta=\{\alpha_1,\dotsc,\alpha_s\}\) be the system of simple roots
    that determines~\(A_+\).  If~\(G\) is semi-simple, these roots form a
    basis of~\(A\).  In general, they are linearly independent, so that we can
    extend them to a basis by certain~\(\alpha_j\) for \(s<j\le r\).  Define a
    vector space isomorphism \(\gamma\colon A\to\R^r\) by \(\gamma(\eta)_j =
    \alpha_j(\eta)\) for \(j=1,\dotsc,r\).  This identifies~\(A_+\) with the
    set of \((x_j)\in\R^r\) with \(x_j\ge0\) for \(1\le j\le s\).
    Equation~\eqref{eq:closure_twopoint} identifies
    \(B\defeq\clos(\{0,\xi\})\) with
    \[
    \gamma(B) = \{(x_j)\in\R^r\mid 0\le x_j\le \gamma(\xi)_j,\ j=1,\dotsc,s\}.
    \]
    
    We may assume \(\alpha_j\in X_*(S)\otimes\Q\subseteq A\), so that
    \(\gamma(\Lambda)\subseteq\Q^r\).  Replacing~\(\gamma\) by
    \(n^{-1}\gamma\) for some \(n\in\N^*\), we can achieve
    \(\Z^r\subseteq\gamma(\Lambda)\).  Hence if \((x_j)\in\gamma(B)\), then
    the truncated vector \(\lfloor (x_j)\rfloor\defeq (\lfloor x_j\rfloor)\)
    belongs to \(\gamma(B\cap\Lambda)\).  It satisfies \(\abs{x_j-\lfloor
      x_j\rfloor}<1\) for all~\(j\).  Since the norm
    \(\norm{\gamma(\eta)}_\infty\) is equivalent to the Euclidean norm
    on~\(A\), we have \(d(\gamma^{-1} \lfloor\gamma(\eta)\rfloor,\eta)<R\) for
    all \(\eta\in B\) for some \(R>0\).
  \end{proof}
  
  The building \(\BT\) can be defined as the quotient of \(G\times A\) by a
  certain equivalence relation.  We may view \(gA\) for \(g\in G\) as a
  subspace of \(\BT\); these are the \emph{apartments} of \(\BT\).  We now
  choose~\(\xi_0\) to be the origin in \(A\subseteq\BT\).  Recall that
  \(H\subseteq G\) denotes the stabiliser of~\(\xi_0\).  We have the
  \emph{Cartan decomposition} \(G=H\Lambda_+H\)
  by~\cite{Tits:Corvallis}*{(3.3.3)}, so that \(G\xi_0=H\Lambda_+\xi_0\).  Let
  \(G_\circ\subseteq G\) be the connected component of the identity (as an
  algebraic variety) and let \(H_\circ\defeq H\cap G_\circ\).  These are open
  normal subgroups of finite index in \(G\) and~\(H\), respectively, and
  \(\BT(G)\) is isomorphic to \(\BT(G_\circ)\) equipped with a canonical
  action of~\(G\).

  Choose \(\Lambda'_+\subseteq\Lambda_+\) to contain one representative for
  each \(H\)\nbd{}orbit in \(G\xi_0\).  Fix \(\xi\in\Lambda'_+\) and
  \(k\in\N\).  We further decompose \(H\xi\) as a disjoint union of finitely
  many \(H_\circ\)\nbd{}orbits \(H_\circ h_j\xi=h_j H_\circ\xi\) for suitable
  \(h_1,\dotsc,h_N\in H\).  We let \(p_k'(\xi)\) be some point in
  \(\Lambda\cap\clos(\{0,\xi\})\subseteq A\subseteq\BT\) that has minimal
  distance from \(p_k(\xi)\).

  \begin{proposition}
    \label{pro:stabilisers}
    Let \(\Omega\subseteq A\subseteq\), \(\Omega\neq\emptyset\).  If \(g\in
    G_\circ\) satisfies \(gx=x\) for all \(x\in\Omega\), then \(gx=x\) for all
    \(x\in\clos(\Omega)\).  (This may fail if we allow \(g\in G\).)
  \end{proposition}

  \begin{proof}
    The proof requires some facts about stabilisers of points in \(\BT\),
    which are conveniently summarised in \cite{Schneider-Stuhler}*{Section
      I.1}.  The subgroups \(P_\Omega\subseteq G_\circ\) defined there
    manifestly satisfy \(P_\Omega=P_{\clos(\Omega)}\).  This implies our claim
    because
    \[
    P_\Omega = \{g\in G_\circ\mid gx=x \ \forall x\in\Omega\}.\qedhere
    \]
  \end{proof}

  Proposition~\ref{pro:stabilisers} allows us to define \(p_k'(h_j h_\circ\xi)
  \defeq h_j h_\circ p_k'(\xi)\) for all \(h_\circ\in H_\circ\).  Letting
  \(\xi,j\) vary, we get a map \(p_k'\colon G\xi_0\to G\xi_0\).  It is
  \(H_\circ\)\nbd{}equivariant because~\(H_\circ\) is normal in~\(H\).  Since
  \(p_k(\xi)\in\clos(\{0,\xi\})\), Lemma~\ref{lem:approx_good} yields \(R>0\)
  such that \(d\bigl(p_k(\xi),p_k'(\xi)\bigr)\le R\) for all \(k\in\N\),
  \(\xi\in\Lambda'_+\).  The same holds for \(\xi\in G\xi_0\) because~\(p_k\)
  is \(H\)\nbd{}equivariant and~\(G\) acts isometrically on \(\BT\).
  Moreover, \(p_0'(\xi)=\xi_0\) and \(p_k'(\xi)=\xi\) for \(k\ge
  d(\xi,\xi_0)\).  Thus \((p_k')\) is a combing of \(G/H\) of linear growth.
  
  \begin{lemma}
    \label{lem:BT_combing_smooth}
    The combing \((p_k')\) is smooth.
  \end{lemma}

  \begin{proof}
    There is a decreasing sequence \((U_n)_{n\in\N}\) in \(\COM(G)\) such that
    each~\(U_n\) is normal in~\(H\) and can be written as \(U_n^+ \cdot U_n^-=
    U_n^-\cdot U_n^+\) with
    \begin{equation}
      \label{eq:congruence_subgroup}
      \lambda U_n^+\subseteq U_n^+\lambda,\qquad
      U_n^-\lambda\subseteq \lambda U_n^-
    \end{equation}
    for all \(\lambda\in\Lambda_+\) (see~\cite{Schneider-Stuhler}*{Section
      I.2}).  Let \(a,b\in G\) satisfy \(p_k'(abH)=aH\).  Write
    \(ab=h_1\lambda h_2\) with \(h_1,h_2\in H\), \(\lambda\in\Lambda_+'\),
    and~\(h_1\) chosen as carefully as above if~\(G\) is disconnected.  Then
    \(a=h_1p_k'(\lambda) h_3\) with the same~\(h_1\) and some \(h_3\in H\).
    Hence \(b=h_3^{-1} p_k'(\lambda)^{-1}\lambda h_2\).  Since
    \(p_k'(\lambda)\in\clos(\{0,\lambda\})\), equation
    \eqref{eq:closure_twopoint} yields \(p_k'(\lambda)\in\Lambda_+\) and
    \(p_k'(\lambda)^{-1}\lambda\in\Lambda_+\).
    Using~\eqref{eq:congruence_subgroup} and that~\(U_n\) is normal in~\(H\),
    we get
    \begin{multline*}
      aU_n b
      = h_1 p_k'(\lambda) U_n^+ U_n^- p_k'(\lambda)^{-1}\lambda h_2
      \subseteq h_1 U_n^+ p_k'(\lambda) p_k'(\lambda)^{-1}\lambda U_n^- h_2
      \\ \subseteq h_1 U_n\lambda U_n h_2
      = U_n h_1\lambda h_2 U_n
      = U_n ab U_n.
    \end{multline*}
    Thus the combing \((p_k')\) is smooth in the sense of
    Definition~\ref{def:smooth_combing}.
  \end{proof}

  Finally, the compatibility condition of
  Definition~\ref{def:compatible_combing} is easy to check using the explicit
  description of \(\Sch(G^2)\) in Lemma~\ref{lem:Sch_tensor}.  In order to
  cover also the Schwartz algebras for discrete groups, which are defined by
  \(\ell_1\)\nbd{}estimates, we define spaces
  \[
  L_p^\sigma (G) \defeq
  \{f\colon G\to\C \mid f\cdot\sigma^k\in L_p(G) \ \forall k\in\N \}
  \]
  for all \(1\le p<\infty\), and equip them with the evident bornology: a
  subset \(T\subseteq L_p^\sigma(G)\) is bounded if and only if for any
  \(k\in\N\) there is~\(C_k\) such that \(\norm{f\cdot\sigma^k}_{L_p(G)}\le
  C_k\) for all \(f\in T\).  Let \(L_p^\sigma(G\biinv U)\) be the subspace of
  \(U\)\nbd{}bi-invariant functions in \(L_p^\sigma(G)\).

  For \(p=2\), this agrees with the previous definition, so that
  Lemma~\ref{lem:Sch_tensor} yields
  \[
  \Sch(G) = \varinjlim L_2^\sigma(G\biinv U),
  \qquad
  \Sch(G^2) = \varinjlim L_2^\sigma(G^2\biinv U^2).
  \]

  \begin{lemma}
    \label{lem:third_condition}
    Let \((p_k)\) be a combing of polynomial growth on \(G/H\).  Then the
    sequence of operators \((R_k)\) used in
    Definition~\ref{def:compatible_combing} is uniformly bounded as operators
    \(L_p^\sigma(G)\to L_p^\sigma(G^2)\).  Here we use the scale
    \(\sigma(a,b)\defeq \sigma(a)\sigma(b)\) on~\(G^2\).
  \end{lemma}

  \begin{proof}
    The operator \(W\phi(x,y)\defeq \phi(x,x^{-1}y)\) is an isometry of
    \(L_p(G^2)\).  It is also sufficiently compatible with the scale
    on~\(G^2\) for~\(W\) and its inverse to be bounded linear operators on
    \(L_p^\sigma(G^2)\).  We have \(R_k\phi(a,b)=\phi(ab) 1_{p_k(abH)}(a)\),
    where \(1_{p_k(abH)}\) denotes the characteristic function of
    \(p_k(abH)\subseteq G\).  Hence
    \[
    WR_k\phi(x,y)=\phi(y) \cdot 1_{p_k(yH)}(x).
    \]
    Since the combing \((p_k)\) has polynomial growth,
    \(\sigma\bigl(p_k(yH)\bigr)\) is controlled by a polynomial in
    \(\sigma(y)\).  The boundedness of \(W\circ R_k\) is now immediate because
    all cosets \(xH\) have volume~\(1\).  This implies the boundedness
    of~\(R_k\).
  \end{proof}
  
  Since the combing \((p_k')\) is smooth, for any \(U\in\COM(G)\) there exists
  \(V\in\COM(G)\) such that~\(R_k\) maps \(\Tc(G\biinv U)\) into
  \(\Tc(G^2\biinv V^2)\).  Together with Lemma~\ref{lem:third_condition} for
  \(p=2\), this yields that the combing \((p_k')\) is compatible with
  \(\Sch(G)\) in the sense of Definition~\ref{def:compatible_combing}.  We
  have now verified all the hypotheses of Theorem~\ref{the:combing_contracts}.
  Thus \(\Sch(G)\) is isocohomological.
\end{proof}

\subsection{Decomposition with respect to the centre of~\(G\)}
\label{sec:centre_decomposition}

As before, we let~\(G\) be a reductive \(p\)\nbd{}adic group.  Let
\(C(G)\subseteq G\) be the connected centre of~\(G\) and let \(\chi\colon
C(G)\to\U(1)\) be a unitary character on~\(C(G)\).  Let \(\Mod_\chi(G)\) be
the full subcategory of \(\Mod(G)\) whose objects are the representations
\(\pi\colon G\to\Aut(V)\) that satisfy \(\pi(z)= \chi(z)\ID_V\) for all \(z\in
C(G)\).  Let
\[
\Mod_\chi\bigl(\Sch(G)\bigr) \defeq \Mod\bigl(\Sch(G)\bigr)\cap \Mod_\chi(G)
\]
be the subcategory of tempered representations in \(\Mod_\chi(G)\).  The class
of extensions with a bounded linear section turns \(\Mod_\chi(G)\) and
\(\Mod_\chi\bigl(\Sch(G)\bigr)\) into exact categories, so that we can form
the derived categories \(\Der_\chi(G)\) and \(\Der_\chi\bigl(\Sch(G)\bigr)\).
Let \(\Gss\defeq G/C(G)\), this is again a reductive \(p\)\nbd{}adic group.
If \(\chi=1\), then \(\Mod_\chi(G)=\Mod(\Gss)\) and
\(\Mod_\chi\bigl(\Sch(G)\bigr)=\Mod\bigl(\Sch(\Gss)\bigr)\).  In general,
there are quasi-unital bornological algebras \(\Hecke_\chi(G)\) and
\(\Sch_\chi(G)\) such that
\[
\Mod_\chi(G) \cong \Mod\bigl(\Hecke_\chi(G)\bigr),
\qquad
\Mod_\chi\bigl(\Sch(G)\bigr) \cong \Mod\bigl(\Sch_\chi(G)\bigr).
\]
We briefly recall their well-known definitions.  A \(C(G)\)\nbd{}invariant
subset of~\(G\) is called \emph{\(C(G)\)\nbd{}compact} if its image
in~\(\Gss\) is compact.  Let \(\Hecke_\chi(G)\) be the space of locally
constant functions \(f\colon G\to\C\) with \(C(G)\)\nbd{}compact support such
that \(f(z^{-1}g)=\chi(z)f(g)\) for all \(g\in G\), \(z\in C(G)\).  If
\(f_1,f_2\in\Hecke_\chi(G)\), then the function \(h\mapsto
f_1(h)f_2(h^{-1}g)\) is \(C(G)\)\nbd{}invariant, so that
\begin{displaymath}
  f_1*f_2(g)\defeq \int_{\Gss} f_1(h)f_2(h^{-1}g)\,dh
\end{displaymath}
makes sense; here~\(dh\) denotes the Haar measure on~\(\Gss\).  This turns
\(\Hecke_\chi(G)\) into an algebra, which we equip with the fine bornology.

Let \(L_2^\sigma(G)_\chi\) be the space of functions \(f\colon G\to\C\) that
satisfy \(f(z^{-1}g)=\chi(z)f(g)\) for all \(g\in G\), \(z\in C(G)\), and such
that \(zg\mapsto \abs{f(g)}\) is an element of \(L_2^\sigma(\Gss)\).  Let
\[
\Sch_\chi(G)\defeq \varinjlim L_2^\sigma(G\biinv U)_\chi
\]
where~\(U\) runs through the set of compact open subgroups with \(U\cap
C(G)\subseteq \ker\chi\).  The same estimates as for \(\Sch(\Gss)\) show that
the convolution on \(\Hecke_\chi(G)\) extends to a bounded multiplication on
\(\Sch_\chi(G)\).

Consider the map
\[
\rho\colon \Hecke(G)\to\Hecke_\chi(G),
\qquad
\rho f(g)\defeq \int_{C(G)} \chi(z) f(zg)\,dz.
\]
For appropriately normalised Haar measures, this is a surjective, bounded
algebra homomorphism; that is, \(\Hecke_\chi(G)\) is a quotient algebra of
\(\Hecke(G)\).  Using~\(\rho\), we can pull back
\(\Hecke_\chi(G)\)\brd{}modules to \(\Hecke(G)\)\brd{}modules.  This
construction maps essential modules again to essential modules (\(\rho\) is a
proper morphism in the notation of~\cite{Meyer:Embed_derived}).  Thus we have
got a functor \(\rho^*\colon
\Mod\bigl(\Hecke_\chi(G)\bigr)\to\Mod\bigl(\Hecke(G)\bigr)\).  Since~\(\rho\)
is surjective, \(\rho^*\) is fully faithful.  Thus
\(\Mod\bigl(\Hecke_\chi(G)\bigr)\) becomes a full subcategory of \(\Mod(G)\).
It is easy to identify this subcategory with \(\Mod_\chi(G)\).  If
\((V,\pi)\in\Mod_\chi(G)\), then~\(V\) becomes an essential
\(\Hecke_\chi(G)\)\brd{}module by
\[
\pi(f) \defeq \int_{\Gss} f(g)\pi(g)\,dg.
\]
This is well-defined because \(f(gz)\pi(gz)=f(g)\pi(g)\) for all \(g\in G\),
\(z\in C(G)\).

We can extend~\(\rho\) to a bounded algebra homomorphism \(\rho^\Sch\colon
\Sch(G)\to\Sch_\chi(G)\).  The map~\(\rho^\Sch\) has a bounded linear section
and its kernel is the closure of \(\ker \rho\subseteq \Hecke(G)\).  Therefore,
bounded algebra homomorphisms \(\Sch_\chi(G)\to\Endo(V)\) correspond to
bounded algebra homomorphisms \(\Sch(G)\to\Endo(V)\) whose restriction to
\(\Hecke(G)\) vanishes on \(\ker\rho\).  Equivalently,
\(\Mod\bigl(\Sch_\chi(G)\bigr) \cong \Mod\bigl(\Sch(G)\bigr)\cap
\Mod_\chi(G)\).  Thus \(\Hecke_\chi(G)\) and \(\Sch_\chi(G)\) have the
required properties.

\begin{theorem}
  \label{the:Sch_chi_isocoh}
  The embedding \(\Hecke_\chi(G)\to\Sch_\chi(G)\) is isocohomological.
\end{theorem}

\begin{proof}
  Let \(\Hecke_\chi(G)^\op\) be the opposite algebra of \(\Hecke_\chi(G)\), so
  that \(\Mod(\Hecke_\chi(G)^\op)\) is the category of \emph{right}
  \(\Hecke_\chi(G)\)\brd{}modules.  Since
  \(\Hecke_\chi(G)^\op\cong\Hecke_{\chi^{-1}}(G)\), we have an isomorphism of
  categories \(\Mod(\Hecke_\chi(G)^\op) \cong \Mod_{\chi^{-1}}(G)\).
  
  Equip \(X\in\Mod(\Hecke_\chi(G)^\op)\),
  \(V\in\Mod\bigl(\Hecke_\chi(G)\bigr)\) with the associated representations
  of~\(G\).  We equip \(X\hot V\) with the diagonal representation.  Since
  \(\chi\) and~\(\chi^{-1}\) cancel, \(C(G)\) acts trivially on \(X\hot V\).
  Thus we obtain a bifunctor
  \begin{equation}
    \label{eq:tensor_chi}
    \Mod(\Hecke_\chi(G)^\op)\times\Mod\bigl(\Hecke_\chi(G)\bigr)
    \to \Mod(\Gss), \qquad (X,V)\mapsto X\hot V.
  \end{equation}
  This functor is evidently exact for extensions with a bounded linear
  section.  Moreover, we claim that \(X\hot V\) is projective if \(X\)
  or~\(V\) are projective.  It suffices to treat the case where~\(X\) is
  projective.  We may even assume that~\(X\) is a free essential module
  \(X_0\hot\Hecke_\chi(G)\).  The diagonal representation on \(X_0\hot
  \Hecke_\chi(G)\hot Y\) is isomorphic to the regular representation
  \(\rho_g\otimes 1\otimes 1\) on \(\Hecke(\Gss)\hot X_0\hot Y\).  The
  intertwining operator is given by \(\Phi (x\otimes f\otimes y)(g)\defeq f(g)
  x\otimes gy\) for all \(g\in G\); this function only depends on the class
  of~\(g\) in~\(\Gss\).
  
  Let \(X\in\Mod(\Hecke_\chi(G)^\op)\),
  \(V\in\Mod\bigl(\Hecke_\chi(G)\bigr)\).  Then \(X\hot_{\Hecke_\chi(G)} V\)
  is defined as the quotient of \(X\hot V\) by the closed linear span of
  \(x*f\otimes v-x\otimes f*v\) for \(x\in X\), \(f\in\Hecke_\chi(G)\), \(v\in
  V\).  Since \(\rho\colon \Hecke(G)\to\Hecke_\chi(G)\) is surjective, this is
  the same as \(X\hot_{\Hecke(G)} V\), which we have identified with \(X\hot_G
  V\) in Section~\ref{sec:homological_Hecke}.  Thus
  \[
  X\hot_{\Hecke_\chi(G)} V \cong \C \hot_\Gss (X\hot V).
  \]
  The same assertion holds for the total derived functors because the
  bifunctor in~\eqref{eq:tensor_chi} is exact and preserves projectives.
  Especially, we get
  \[
  \Sch_\chi(G)\Lhot_{\Hecke_\chi(G)} \Sch_\chi(G)
  \cong \C \Lhot_\Gss \bigl(\Sch_\chi(G)\hot\Sch_\chi(G)\bigr).
  \]
  Here \(\Sch_\chi(G)\hot\Sch_\chi(G)\) is equipped with the inner conjugation
  action of~\(\Gss\).  We identify \(\Sch_\chi(G)\hot\Sch_\chi(G) \cong
  \Sch_{\chi\times\chi}(G\times G)\) as in Lemma~\ref{lem:Sch_tensor}.  By
  \cite{Meyer:Embed_derived}*{Theorem 35.2}, the embedding
  \(\Hecke_\chi(G)\to\Sch_\chi(G)\) is isocohomological if and only if
  \(\Sch_\chi(G)\Lhot_{\Hecke_\chi(G)} \Sch_\chi(G)\) is a resolution of
  \(\Sch_\chi(G)\).  Thus the assertion that we have to prove is equivalent to
  \[
  \C \Lhot_\Gss \Sch_{\chi\times\chi}(G\times G)
  \cong \Sch_\chi(G).
  \]
  We already know \(\C\Lhot_\Gss \Sch(\Gss\times\Gss)\cong\Sch(\Gss)\) because
  \(\Sch(\Gss)\) is isocohomological (Theorem~\ref{the:Sch_isocohomological})
  and this condition is equivalent to \(\Sch(\Gss)\) being isocohomological.

  Now we choose a continuous section \(s\colon \Gss\to G\); this is possible
  because~\(G\) is totally disconnected.  It yields bornological isomorphisms
  \begin{alignat*}{2}
    \Psi'&\colon \Sch_\chi(G)\to\Sch(\Gss),
    &\qquad \Psi' f(g)&\defeq f\circ s(g),
    \\
    \Psi &\colon \Sch_{\chi\times\chi}(G\times G) \to \Sch(\Gss\times\Gss),
    &\qquad \Psi f(g,h) &\defeq f\bigl(s(g),s(g)^{-1} s(gh)\bigr).
  \end{alignat*}
  The isomorphism~\(\Psi\) intertwines the inner conjugation actions
  of~\(\Gss\) on \(\Sch_{\chi\times\chi}(G\times G)\) and
  \(\Sch(\Gss\times\Gss)\).  Thus we get isomorphisms
  \begin{equation}
    \label{eq:Sch_chi_iso}
    \C\Lhot_\Gss \Sch_{\chi\times\chi}(G\times G)
    \cong \C\Lhot_\Gss \Sch(\Gss\times\Gss)
    \cong \Sch(\Gss)
    \cong \Sch_\chi(G).
  \end{equation}
  It is easy to see that the composite isomorphism is induced by the
  convolution map in \(\Sch_\chi(G)\).
\end{proof}

\section{Applications to representation theory}
\label{sec:applications}

Let~\(G\) be a reductive \(p\)\nbd{}adic group, let \(C(G)\) be its centre,
and let \(\Gss\defeq G/C(G)\).  So far we have used very large projective
\(\Hecke(G)\)\brd{}module resolutions, which offer great flexibility for
writing down contracting homotopies.  Now we consider much smaller projective
resolutions, which are useful for explicit calculations.  We write
\(\Mod_{(\chi)}(G)\) if it makes no difference whether we work in \(\Mod(G)\)
or \(\Mod_\chi(G)\) for some character \(\chi\colon C(G)\to\U(1)\).
Similarly, we write \(\Hecke_{(\chi)}(G)\) and \(\Sch_{(\chi)}(G)\).  The
actual applications of our main theorem are contained in Sections
\ref{sec:coh_dim}, \ref{sec:dim_si} and~\ref{sec:vanishing}.  The other
subsections contain small variations on known results.  Our presentation
differs somewhat from the accounts in \cites{Schneider-Stuhler,
  Vigneras:Dimension} because we want to exhibit connections with
\(\K\)\nbd{}theory and assembly maps.

\subsection{Cohomological dimension}
\label{sec:coh_dim}

Let \(\rk G=\dim \BT(G)\) be the \emph{rank} of~\(G\).

\begin{theorem}
  \label{the:cohom_dim}
  The cohomological dimensions of the exact categories \(\Mod(G)\) and
  \(\Mod\bigl(\Sch(G)\bigr)\) are (at most) \(\rk G\); that is, any object has
  a projective resolution of length \(\rk G\).  Similarly, the cohomological
  dimensions of \(\Mod_\chi(G)\) and \(\Mod_\chi\bigl(\Sch(G)\bigr)\) for a
  character \(\chi\colon C(G)\to\U(1)\) are at most \(\rk\Gss\).
\end{theorem}

\begin{proof}
  The assertions are well-known for \(\Mod(G)\) and \(\Mod_\chi(G)\).  For the
  proof, equip \(\BT=\BT(G)\) with a CW-complex structure for which~\(G\) acts
  by cellular maps.  Then the cellular chain complex \(C_\bullet(\BT)\) is a
  projective \(\Hecke(G)\)\brd{}module resolution of the trivial
  representation of length \(\rk G\).  The chain complex \(C_\bullet(\BT)\hot
  V\) with the diagonal representation of~\(G\) is a projective
  \(\Hecke(G)\)\brd{}module resolution of~\(V\) for arbitrary \(V\in\Mod(G)\).
  If \(V\in\Mod_\chi(G)\), we use the building \(\BT(\Gss)\) instead;
  \(C_\bullet\bigl(\BT(\Gss)\bigr)\hot V\) is a projective resolution of~\(V\)
  in \(\Mod_\chi(G)\).

  What is new is that we get the same assertions for
  \(\Mod\bigl(\Sch(G)\bigr)\) and \(\Mod\bigl(\Sch_\chi(G)\bigr)\).  Since the
  argument is the same in both cases, we only write it down for \(\Sch(G)\).
  Let \(V\in\Mod\bigl(\Sch(G)\bigr)\) and let \(P_\bullet\to V\) be a
  projective resolution in \(\Mod(G)\) of length \(\rk G\).  Then
  \(\Sch(G)\hot_G P_\bullet\) has the homotopy type of
  \begin{equation}
    \label{eq:Sch_Lhot_V}
    \Sch(G)\Lhot_G V\cong \Sch(G)\Lhot_{\Sch(G)} V \cong V    
  \end{equation}
  because \(\Sch(G)\) is isocohomological
  (Theorem~\ref{the:Sch_isocohomological}); here we use one of the equivalent
  characterisations of isocohomological embeddings listed in
  Section~\ref{sec:isocoh}.  Equation~\eqref{eq:Sch_Lhot_V} means that
  \(\Sch(G)\hot_G P_\bullet\) is a resolution of~\(V\).  This resolution is
  projective and has length \(\rk G\).
\end{proof}

Conversely, there is \(V\in\Mod(G)\) with \(\Ext^{\rk G}_G(V,V)\neq0\).  Hence
the cohomological dimension of \(\Mod(G)\) is equal to \(\rk G\).  We can even
take~\(V\) tempered and irreducible.  Hence \(\Ext^{\rk
  G}_{\Sch(G)}(V,V)\neq0\) as well because \(\Sch(G)\) is isocohomological.
Thus \(\Mod\bigl(\Sch(G)\bigr)\) also has cohomological dimension equal to
\(\rk G\).  Similarly, the cohomological dimension of \(\Mod_\chi(G)\) and
\(\Mod_\chi\bigl(\Sch(G)\bigr)\) is equal to \(\rk \Gss\).

\subsection{Finite projective resolutions}
\label{the:finite_resolution}

We use a result of Joseph Bernstein to attach an Euler characteristic
\(\Eul(V)\) in \(\K_0\bigl(\Hecke_{(\chi)}(G)\bigr)\) to a finitely generated
representation \(V\in\Mod_{(\chi)}(G)\).

\begin{definition}
  \label{def:finitely_generated}
  A smooth representation~\(V\) is called \emph{finitely generated} if there
  exist finitely many elements \(v_1,\dotsc,v_n\) such that the map
  \[
  \Hecke(G)^n\to V,
  \qquad (f_1,\dotsc,f_n) \mapsto \sum_{j=1}^n f_j*v_j
  \]
  is a bornological quotient map.
\end{definition}

An admissible representation is finitely generated if and only if it has
\emph{finite length}, that is, it has a Jordan-Hölder series of finite length.

Since \(\Hecke(G)^n\) carries the fine bornology, the same is true for its
quotients.  Hence a finitely generated representation necessarily belongs to
\(\Mod_\alg(G)\).  In the situation of Definition
\ref{def:finitely_generated}, there exists \(U\in\COM(G)\) fixing~\(v_j\) for
all \(j\in\{1,\dotsc,n\}\).  Thus we get a bornological quotient map
\(\Hecke(G/U)^n\prto V\).  Conversely, \(\Hecke(G/U)^n\) is finitely generated
and projective.  Thus a smooth representation is finitely generated if and
only if it is a quotient of \(\Hecke(G/U)^n\) for some \(U\in\COM(G)\).  If
\(V\in\Mod_\chi(G)\), then we may replace \(\Hecke(G/U)^n\) by
\(\Hecke_\chi(G/U)^n\) for some \(U\in\COM(G)\) with \(\chi|_{U\cap C(G)}=1\).

\begin{definition}
  \label{def:FP}
  An object of \(\Mod_{(\chi)}(G)\) has type \emph{(FP)} if it admits a
  resolution of finite length by finitely generated projective objects of
  \(\Mod_{(\chi)}(G)\).  Such a resolution is called a \emph{finite projective
    resolution}.
\end{definition}

\begin{theorem}[Joseph Bernstein]
  \label{the:Bernstein}
  An object of \(\Mod_{(\chi)}(G)\) has type (FP) if and only if it is
  finitely generated.
\end{theorem}

\begin{proof}
  It is trivial that representations of type (FP) are finitely generated.
  Conversely, if~\(V\) is finitely generated, then~\(V\) is a quotient of a
  finitely generated projective representation, say, \(\partial_0\colon
  \Hecke_{(\chi)}(G/U)^n\prto V\).  By \cite{Bernstein:centre}*{Remark 3.12},
  subrepresentations of finitely generated representations are again finitely
  generated.  Especially, \(\ker \partial_0\) is finitely generated.  By
  induction, we get a resolution \((P_n,\partial_n)\) of~\(V\) by finitely
  generated projective objects.  By Theorem~\ref{the:cohom_dim}, the kernel of
  \(\partial_n\colon P_n\to P_{n-1}\) is projective for sufficiently
  large~\(n\).  Hence
  \[
  0 \to \ker \partial_n \to P_n \to \dotso \to P_0\to V
  \]
  is a finite projective resolution.
\end{proof}

The \emph{algebraic \(\K\)\nbd{}theory} \(\K_0\bigl(\Hecke_{(\chi)}(G)\bigr)\)
is the Grothendieck group of the monoid of finitely generated projective
\(\Hecke_{(\chi)}(G)\)\brd{}modules.  This is so because
\(\Hecke_{(\chi)}(G)\) is a union of unital subalgebras.

\begin{definition}
  \label{def:Euler_characteristic}
  Let \(V\in\Mod_{(\chi)}(G)\) be finitely generated.  Then~\(V\) is of type
  (FP) by Bernstein's Theorem~\ref{the:Bernstein}.  Choose a finite projective
  resolution
  \[
  0 \to P_n \to \dotso \to P_0 \to V \to 0.
  \]
  The \emph{Euler characteristic} of~\(V\) is defined by
  \[
  \Eul(V) \defeq \sum_{j=0}^n (-1)^j [P_j] \in
  \K_0\bigl(\Hecke_{(\chi)}(G)\bigr).
  \]
\end{definition}

We check that this does not depend on the resolution (see also
\cite{Rosenberg:Algebraic_K}*{Section 1.7}).  Define the Euler characteristic
\(\Eul(P_\bullet)\) for finite projective complexes in the obvious fashion.
Let \(P_\bullet\) and~\(P'_\bullet\) be two finite projective resolutions
of~\(V\).  The identity map on~\(V\) lifts to a chain homotopy equivalence
\(f\colon P_\bullet\to P'_\bullet\).  Hence the mapping cone~\(C_f\) of~\(f\)
is contractible.  The Euler characteristic vanishes for contractible
complexes.  Hence \(\Eul(C_f)=0\).  This is equivalent to
\(\Eul(P_\bullet)=\Eul(P'_\bullet)\).

\begin{definition}
  \label{def:universal_trace}
  Let
  \[
  \HH_0\bigl(\Hecke_{(\chi)}(G)\bigr)
  \defeq \Hecke_{(\chi)}(G)/[\Hecke_{(\chi)}(G),\Hecke_{(\chi)}(G)].
  \]
  The \emph{universal trace} is a map
  \[
  \tr_\univ\colon \K_0\bigl(\Hecke_{(\chi)}(G)\bigr) \to
  \HH_0\bigl(\Hecke_{(\chi)}(G)\bigr).
  \]
  If \((p_{ij}) \in M_n\bigl(\Hecke_{(\chi)}(G)\bigr)\) is an idempotent with
  \(\Hecke_{(\chi)}(G)^n\cdot (p_{ij}) \cong V\), then we have
  \(\tr_\univ[V]=\bigl[\sum p_{ii}\bigr]\).
\end{definition}

The above definitions are inspired by constructions of Hyman Bass
in~\cite{Bass:Euler}, where \(\tr_\univ \Eul(V)\) is constructed for modules
of type (FP) over unital algebras.

\subsection{Traces from admissible representations}
\label{sec:traces_admissible}

Let \(W\in\Mod_{(\chi)}(G)\) be an admissible representation.  Its integrated
form is an algebra homomorphism~\(\rho\) from \(\Hecke_{(\chi)}(G)\) to the
algebra \(\Endo_\fin(W)\defeq \tilde{W}\otimes W\) of smooth finite rank
operators on~\(W\).  Here~\(\tilde{W}\) denotes the contragradient
representation and \(\Endo_\fin(W)\) carries the fine bornology.
Composing~\(\rho\) with the standard trace on \(\Endo_\fin(W)\), we get a
trace \(\tr_W\colon \HH_0\bigl(\Hecke_{(\chi)}(G)\bigr)\to\C\) and a
functional
\[
\tau_W\colon \K_0\bigl(\Hecke_{(\chi)}(G)\bigr)\to\Z.
\]
The following computation of~\(\tau_W\) is a variant of
\cite{Bass:Euler}*{Proposition 4.2}.

\begin{proposition}
  \label{pro:tau_WV}
  Let \(V,W\in\Mod_{(\chi)}(G)\), let~\(V\) be finitely generated projective
  and let~\(W\) be admissible.  Then \(\tau_W[V] = \dim \Hom_G(V,W)\) and
  \(\Hom_G(V,W)\) is finite-dimensional.
\end{proposition}

\begin{proof}
  The functoriality of \(\K_0\) for the homomorphism
  \(\Hecke_{(\chi)}(G)\to\Endo_\fin(W)\) maps \([V]\) to the class of the
  finitely generated projective module
  \[
  V'\defeq \Endo_\fin(W)\hot_{\Hecke_{(\chi)}(G)} V
  \cong W\otimes (\tilde{W}\otimes_G V)
  \cong W^{\dim\tilde{W}\otimes_G V}
  \]
  over \(\Endo_\fin(W)\).  Thus \(\tau_W[V]=\dim \tilde{W}\otimes_G V\). By
  adjoint associativity,
  \[
  \Hom(\tilde{W}\otimes_G V,\C) \cong \Hom_G\bigl(V,\Hom(\tilde{W},\C)\bigr)
  \cong \Hom_G\bigl(V,\tilde{\tilde{W}}\bigr).
  \]
  We have \(W\cong\tilde{\tilde{W}}\) because~\(W\) is admissible.  Thus
  \(\tau_W[V]= \dim \Hom_G(V,W)\).
\end{proof}

Let \(V,W\in\Mod_{(\chi)}(G)\), let~\(W\) be admissible and let~\(V\) be
finitely generated.  By Bernstein's Theorem~\ref{the:Bernstein}, there is a
finite projective resolution \(P_\bullet\to V\) in \(\Mod_{(\chi)}(G)\).  By
Proposition~\ref{pro:tau_WV}, \(\Hom_G(P_\bullet,W)\) is a chain complex of
finite-dimensional vector spaces.  Hence its homology
\(\Ext_{\Hecke_{(\chi)}(G)}^n(V,W)\) is finite-dimensional as well and
\begin{multline}
  \label{eq:two_sums}
  \sum_{n=0}^\infty (-1)^n \dim \Ext_{\Hecke_{(\chi)}(G)}^n(V,W)
  = \sum_{n=0}^\infty (-1)^n \dim \Hom_G(P_n,W)
  \\ = \sum_{n=0}^\infty (-1)^n \tau_W[P_n]
  = \tau_W\bigl(\Eul(V)\bigr).
\end{multline}
We call this the \emph{Euler-Poincaré characteristic} \(\EP_{(\chi)}(V,W)\) of
\(V\) and~\(W\) (compare \cite{Schneider-Stuhler}*{page 135}).

\subsection{Formal dimensions}
\label{sec:formal_dim}

Evaluation at \(1\in G\) is a trace on \(\Hecke_{(\chi)}(G)\), that is, a
linear functional \(\tau_1\colon \HH_0\bigl(\Hecke_{(\chi)}(G)\bigr)\to\C\).
The functional
\[
\dim\defeq \tau_1\circ\tr_\univ\colon \K_0\bigl(\Hecke_{(\chi)}(G)\bigr)\to\C
\]
computes the \emph{formal dimension} for finitely generated projective
\(\Hecke_{(\chi)}(G)\)\brd{}modules.  Recall that an irreducible
representation in \(\Mod_\chi(G)\) is projective if and only if it is
supercuspidal.  Unless \(C(G)\) is compact, \(\Mod(G)\) has no irreducible
projective objects.

We can also define the formal dimension for representations that are
square-integrable (see \cite{Silberger:Reductive_padic}).  Let
\((V,\pi)\in\Mod_\chi(G)\) be irreducible and square-integrable (or, more
precisely, square-integrable \emph{modulo the centre \(C(G)\)}).  Since
irreducible representations are admissible, \(V\) carries the fine bornology.
Moreover, square-integrable representations are tempered.  Thus the integrated
form of~\(\pi\) extends to a bounded homomorphism \(\pi\colon
\Sch_\chi(G)\to\Endo_\fin(V)\).  Since~\(V\) is irreducible, this homomorphism
is surjective.  The crucial property of irreducible square-integrable
representation is that there is an ideal \(I\subseteq\Sch_\chi(G)\) such that
\(I\oplus\ker\pi\cong\Sch_\chi(G)\).  Thus \(\pi|_I\colon I\to\Endo_\fin(V)\)
is an algebra isomorphism.  It is necessarily a bornological isomorphism
because it is bounded and \(\Endo_\fin(V)\) carries the fine bornology.

\begin{proposition}
  \label{pro:si_projective}
  Let \(V\in\Mod_\chi(G)\) be irreducible and square-integrable.  Then~\(V\)
  is both projective and injective as an object of
  \(\Mod_\chi\bigl(\Sch(G)\bigr)\).
\end{proposition}

\begin{proof}
  The direct-sum decomposition \(\Sch_\chi(G) \cong \ker\pi\oplus
  \Endo_\fin(V)\) gives rise to an equivalence of exact categories
  \[
  \Mod_\chi\bigl(\Sch(G)\bigr)
  \cong \Mod(\ker\pi)\times\Mod\bigl(\Endo_\fin(V)\bigr).
  \]
  The representation~\(V\) belongs to the second factor.  The algebra
  \(\Endo_\fin(V)\) is canonically Morita equivalent to~\(\C\), so that
  \(\Mod\bigl(\Endo_\fin(V)\bigr)\) and \(\Mod(\C)\) are equivalent exact
  categories.  The easiest way to get this Morita equivalence uses a basis
  in~\(V\) to identify \(\Endo_\fin(V)\cong\bigcup_{n=1}^\infty M_n(\C)\).
  Since any extension in \(\Mod(\C)\) splits, any object of \(\Mod(\C)\) is
  both injective and projective.
\end{proof}

We can also define \(\K_0\bigl(\Sch_{(\chi)}(G)\bigr)\),
\(\HH_0\bigl(\Sch_{(\chi)}(G)\bigr)\), and
\[
\tr_\univ\colon \K_0\bigl(\Sch_{(\chi)}(G)\bigr) \to
\HH_0\bigl(\Sch_{(\chi)}(G)\bigr).
\]
It is irrelevant for the following whether we divide by the linear or closed
linear span of the commutators in the definition of
\(\HH_0\bigl(\Sch_{(\chi)}(G)\bigr)\).  The trace~\(\tau_1\) extends to a
bounded trace \(\tau^\Sch_1\colon\HH_0\bigl(\Sch_{(\chi)}(G)\bigr)\to\C\).
This induces a functional \(\dim^\Sch\colon
\K_0\bigl(\Sch_{(\chi)}(G)\bigr)\to\Z\).  An irreducible square-integrable
representation~\(V\) defines a class
\([V]\in\K_0\bigl(\Sch_{(\chi)}(G)\bigr)\) by
Proposition~\ref{pro:si_projective}; we define its \emph{formal dimension} by
\(\dim^\Sch V \defeq \dim^\Sch[V]\).

The embedding \(\Hecke_{(\chi)}(G)\to\Sch_{(\chi)}(G)\) induces natural maps
\begin{alignat*}{2}
  \iota&\colon \K_0\bigl(\Hecke_{(\chi)}(G)\bigr)\to
  \K_0\bigl(\Sch_{(\chi)}(G)\bigr),
  &\qquad [V]&\mapsto [\Sch_{(\chi)}(V)\hot_{\Hecke_{(\chi)}(G)} V],
  \\
  \iota &\colon \HH_0\bigl(\Hecke_{(\chi)}(G)\bigr)\to
  \HH_0\bigl(\Sch_{(\chi)}(G)\bigr).
  &\qquad [f]&\mapsto [f].
\end{alignat*}
These maps are compatible with the universal traces and satisfy
\(\tau_1^\Sch\circ\iota=\tau_1\) and \({\dim^\Sch}\circ\iota=\dim\).

It is shown in~\cite{Vigneras:Dimension} that \(\Sch_{(\chi)}(G)\) is closed
under holomorphic functional calculus in the \(C^*\)\brd{}algebra
\(C^*_{\mathrm{red},(\chi)}(G)\).  Hence
\[
\K_0\bigl(\Sch_{(\chi)}(G)\bigr) \cong
\K_0\bigl(C^*_{\mathrm{red},(\chi)}(G)\bigr).
\]
It follows also that any finitely generated projective module~\(V\) over
\(\Sch_{(\chi)}(G)\) is the range of a \emph{self-adjoint} idempotent element
in \(M_n\bigl(\Sch_{(\chi)}(G)\bigr)\) for some \(n\in\N\).  Since the
trace~\(\tau_1^\Sch\) is positive, we get \(\dim^\Sch V>0\) unless \(V=0\).

Yet another notion of formal dimension comes from the theory of von Neumann
algebras.  The \emph{(\(\chi\)-twisted) group von Neumann algebra} of~\(G\) is
the closure \(N_{(\chi)}(G)\) of \(\Hecke_{(\chi)}(G)\) or
\(\Sch_{(\chi)}(G)\) in the weak operator topology on \(L_2(G)_{(\chi)}\).  We
may extend~\(\tau_1^\Sch\) to a positive unbounded trace \(\tau_1^N\) on
\(N_{(\chi)}(G)\).  Any normal \(*\)\nbd{}representation~\(\rho\) of
\(N_{(\chi)}(G)\) on a separable Hilbert space is isomorphic to the left
regular representation on the Hilbert space
\(\bigl(L_2(G)_{(\chi)}\barotimes\ell_2(\N)\bigr)\cdot p_\rho\) for some
projection \(p_\rho\in N_{(\chi)}(G)\barotimes B(\ell_2\N)\)
where~\(\barotimes\) denotes spatial tensor products of Hilbert spaces and von
Neumann algebras, respectively; the projection~\(p_\rho\) is unique up to
unitary equivalence.  We define the \emph{formal dimension} \(\dim^N(\rho)\)
to be \(\tau_1^N(p_\rho)\in [0,\infty]\); this does not depend on the choice
of~\(p_\rho\).

By definition, we have \(\dim^N(L_2(G)_{(\chi)}^n\cdot e)=\dim^\Sch[e]\) if
\(e\in M_n\bigl(\Sch_{(\chi)}(G)\bigr)\) is a \emph{self-adjoint} idempotent.
Since \(\Sch_{(\chi)}(G)\) is closed under holomorphic functional calculus in
the reduced group \(C^*\)-algebra \(C^*_{\mathrm{red},(\chi)}(G)\), any
idempotent element of \(\Sch_{(\chi)}(G)\) is similar to a self-adjoint
idempotent.  Hence \(\dim^N(L_2(G)_{(\chi)}^n\cdot e)=\dim^\Sch[e]\) holds for
any idempotent \(e\in M_n\bigl(\Sch_{(\chi)}(G)\bigr)\).  We have
\[
L_2(G)_{(\chi)} \hot_{\Sch_{(\chi)}(G)} \Sch_{(\chi)}(G)^n\cdot e \cong
L_2(G)_{(\chi)}^n \cdot e.
\]
Therefore, if~\(V\) is a finitely generated projective left
\(\Sch_{(\chi)}(G)\)-module~\(V\), then we may view \(L_2(G)_{(\chi)}
\hot_{\Sch_{(\chi)}(G)} V\) as a Hilbert space equipped with a faithful normal
\(*\)\nbd{}representation of \(N_{(\chi)}(G)\); the resulting representation
is uniquely determined up to unitary equivalence because any two self-adjoint
idempotents realising~\(V\) are unitarily equivalent in \(\Sch_{(\chi)}(G)\).
The formal dimensions from the Schwartz algebra and the von Neumann algebra
are compatible in the following sense:
\[
\dim^N(L_2(G)_{(\chi)} \hot_{\Sch_{(\chi)}(G)} V) \cong \dim^\Sch(V).
\]
Thus \(\dim^\Sch V\) only depends on the unitary equivalence class of the
associated unitary representation \(L_2(G)_{(\chi)}\hot_{\Sch_{(\chi)}(G)}
V\).

\subsection{Compactly induced representations}
\label{sec:ci}

Equip \(U\in\COM(G)\) with the restriction of the Haar measure from~\(G\).
The map \(i_U^G\colon \Hecke(U)\to\Hecke(G)\) that extends functions by~\(0\)
outside~\(U\) is an algebra homomorphism.  Hence it induces a map
\[
(i_U^G)_!\colon \Rep(U)\cong \K_0\bigl(\Hecke(U)\bigr)\to\K_0\bigl(\Hecke(G)\bigr),
\qquad
[V]\mapsto [\Hecke(G) \otimes_{\Hecke(U)} V].
\]
This is the standard functoriality of \(\K\)\nbd{}theory.  We denote it by
\((i_U^G)_!\) because this notation is used in~\cite{Meyer:Embed_derived}.  We
call representations of the form \((i_U^G)_!(V)\) \emph{compactly induced}
because \(\Hecke(G)\otimes_{\Hecke(U)} V\cong \cInd_U^G(V)\) (see
\cite{Meyer:Smooth}).

Let \(U,V\in\COM(G)\) and suppose that \(gUg^{-1}\subseteq V\) for some \(g\in
G\).  Then we have \(i_U^G = \gamma_g^{-1}\circ i_V^G\circ
i_{gUg^{-1}}^{V}\circ \gamma_g\), where~\(\gamma_g\) denotes conjugation
by~\(g\).  One checks that~\(\gamma_g\) acts trivially on
\(\K_0\bigl(\Hecke(G)\bigr)\).  Hence \((i_U^G)_!\) is the composite of
\((i_V^G)_!\) and the map \(\Rep(U)\congto\Rep(gUg^{-1})\to \Rep(V)\) that is
associated to the group homomorphism \(U\to V\), \(x\mapsto gxg^{-1}\).  Let
\(\mathsf{Sub}(G)\) be the category whose objects are the compact open
subgroups of~\(G\) and whose morphisms are these special group homomorphisms.
We have exhibited that \(U\mapsto \Rep(U)\) is a module over this category.
The various maps \((i_U^G)_!\) combine to a natural map
\[
\varinjlim_{\mathsf{Sub}(G)} \Rep(U)\to\K_0\bigl(\Hecke(G)\bigr).
\]
We call this map the \emph{assembly map} for \(\K_0\bigl(\Hecke(G)\bigr)\)
because it is a variant of the Farrell-Jones assembly map for discrete groups
(see~\cite{Lueck-Reich:Survey_assembly}*{Conjecture 3.3}), which is in turn
closely related to the Baum-Connes assembly map.

The above definitions carry over to \(\Hecke_\chi(G)\) in a straightforward
fashion.  Let \(\COM\bigl(G;C(G)\bigr)\) be the set of \(C(G)\)\brd{}compact
open subgroups of~\(G\) containing \(C(G)\).  The projection to \(\Gss\)
identifies \(\COM\bigl(G;C(G)\bigr)\) with \(\COM(\Gss)\).  As above, we get
algebra homomorphisms \(i_U^G\colon \Hecke_\chi(U)\to\Hecke_\chi(G)\) for
\(U\in\COM\bigl(G;C(G)\bigr)\).  There is an analogue of the Peter-Weyl
theorem for \(\Hecke_\chi(U)\); that is, \(\Hecke_\chi(U)\) is a direct sum of
matrix algebras.  Therefore, finitely generated projective modules over
\(\Hecke_\chi(U)\) are the same as finite-dimensional representations in
\(\Mod_\chi(U)\).  This justifies defining \(\Rep_\chi(U)\defeq
\K_0\bigl(\Hecke_\chi(U)\bigr)\).  As above, we can factor \((i_U^G)_!\)
through \((i_V^G)_!\) if~\(U\) is subconjugate to~\(V\).  The relevant
category organising these subconjugations is the category
\(\mathsf{Sub}\bigl(G;C(G)\bigr)\) whose set of objects is
\(\COM\bigl(G;C(G)\bigr)\) and whose morphisms are the group homomorphisms
\(U\to V\) of the form \(x\mapsto gxg^{-1}\) for some \(g\in G\).  Thus we get
an assembly map
\begin{equation}
  \label{eq:FJ_assembly}
  \varinjlim_{\mathsf{Sub}\bigl(G;C(G)\bigr)}
  \Rep_\chi(U)\to\K_0\bigl(\Hecke_\chi(G)\bigr).
\end{equation}

Let \(\Gamma(\Gss,dg)\subseteq\R\) be the subgroup generated by
\(\vol(U)^{-1}\) for \(U\in\COM(\Gss)\).  Since \(\vol(U)/\vol(V)\in\N\) for
\(V\subseteq U\), this group is already generated by \(\vol(U)^{-1}\) for
maximal compact subgroups \(U\subseteq\Gss\).  We have
\(\Gamma(\Gss,dg)=\alpha\Z\) for some \(\alpha>0\) because there are only
finitely many maximal compact subgroups and \(\vol(U)/\vol(V)\in\Q\) for all
\(U,V\in\COM(\Gss)\).  The number~\(\alpha\) depends on the choice of the Haar
measure, of course.  We let \(\size(U)\defeq \alpha\vol\bigl(U/C(G)\bigr)\),
so that \(\size(U)^{-1}\in\N\) for all \(U\in\COM\bigl(G;C(G)\bigr)\).

\begin{lemma}
  \label{lem:truniv_cInd}
  Let \(U\in\COM\bigl(G;C(G)\bigr)\), and let \(W\in\Mod_\chi(U)\) be
  finite-dimensional.  Let \(c_W\colon U\to\C\) be the character of~\(W\).
  Then \(\tr_\univ (i_U^G)_! [W] \in \HH_0\bigl(\Hecke_\chi(G)\bigr)\) is
  represented by the function
  \begin{equation}
    \label{eq:induced_representative}
    c_{U,W}^G (g) \defeq
    \begin{cases}
      \alpha \size(U)^{-1} \conj{c_W(g)} & \text{for \(g\in U\),} \\
      0 & \text{for \(g\notin U\).} \\
    \end{cases}
  \end{equation}
  Moreover, \(\dim (i_U^G)_! [W] = \alpha \size(U)^{-1} \dim(W)\).  Thus
  \(\dim x\in\alpha\Z\) for all~\(x\) in the range of the assembly
  map~\eqref{eq:FJ_assembly}.
\end{lemma}

A similar result holds for compactly induced projective objects of
\(\Mod(G)\).

\begin{proof}
  Since the universal trace is compatible with the functoriality of \(\K_0\)
  and \(\HH_0\), the first assertion follows if \(\tr_\univ[W] =
  \vol(U/C(G))^{-1}\conj{c_W(g)}\) in \(\HH_0\bigl(\Hecke_\chi(U)\bigr)\).  We
  briefly recall how this well-known identity is proved.  We may assume
  that~\(W\) is irreducible.  Hence there is an idempotent
  \(p_W\in\Hecke_\chi(U)\) with \(W\cong \Hecke_\chi(U)p_W\).  Thus
  \(\tr_\univ[W]=[p_W]\).  We can compute \(c_W(g)\) for \(g\in G\) as the
  trace of the finite rank operator \(f\mapsto \lambda(g) f*p_W\) on
  \(\Hecke_\chi(U)\).  This operator has the integral kernel \((x,y)\mapsto
  p_W(y^{-1}g^{-1}x)\), so that
  \[
  \conj{c_W(g)} = c_W(g^{-1})
  = \int_{U/C(G)} p_W(x^{-1}gx)\,dx.
  \]
  This implies \([\conj{c_W}] = \int_{U/C(G)} [W]\,dx =
  \vol\bigl(U/C(G)\bigr)[W]\) because conjugation does not change the class in
  \(\HH_0\bigl(\Hecke_\chi(U)\bigr)\).  We get the formula for formal
  dimensions because \(\dim x=\tr_\univ(x)(1)\) and \(c_W(1)=\dim W\).  This
  lies in \(\alpha\Z\) by construction of~\(\alpha\).  By additivity, we get
  \(\dim x \in \alpha\Z\) for all~\(x\) in the range of the assembly
  map~\eqref{eq:FJ_assembly}.
\end{proof}

\subsection{Explicit finite projective resolutions}
\label{sec:explicit_resolution}

Let \(V\in\Mod_\chi(G)\) be of finite length.  Peter Schneider and Ulrich
Stuhler construct an explicit finite projective resolution for such~\(V\)
in~\cite{Schneider-Stuhler}.  We only sketch the construction very briefly.
Let \(\BT(\Gss)\) be the affine Bruhat-Tits building of~\(\Gss\).  One defines
a coefficient system \(\gamma_e(V)\) on \(\BT(\Gss)\), which depends on an
auxiliary parameter \(e\in\N\) (see \cite{Schneider-Stuhler}*{Section II.2});
its value on a facet~\(F\) of \(\BT(\Gss)\) is the\mdash
finite-dimensional\mdash space \(\Fix(U_F^e,V)\) for certain
\(U_F^e\in\COM(G)\).  The cellular chain complex
\(C_\bullet\bigl(\BT(\Gss),\gamma_e(V)\bigr)\) with values in \(\gamma_e(V)\)
is a resolution of~\(V\) for sufficiently large~\(e\)
(\cite{Schneider-Stuhler}*{Theorem II.3.1}).  It is a finite projective
resolution of~\(V\) in \(\Mod_\chi(G)\) because the stabilisers of facets
belong to \(\COM\bigl(G;C(G)\bigr)\) and the set of facets is
\(\Gss\)\nbd{}finite.

\begin{proposition}
  \label{pro:fEP}
  If \(V\in\Mod_\chi(G)\) has finite length, then \(\Eul(V)\) belongs to the
  range of the assembly map~\eqref{eq:FJ_assembly}.  Hence
  \(\dim\bigl(\Eul(V)\bigr)\in\alpha\Z\).
  
  Define the \emph{Euler-Poincaré function} \(f_\EP^V\in\Hecke_\chi(G)\)
  of~\(V\) as in~\cite{Schneider-Stuhler}*{page 135}.  Then \([f_\EP^V] =
  \tr_\univ \Eul(V) \in \HH_0\bigl(\Hecke_\chi(G)\bigr)\).  Thus
  \(\dim\bigl(\Eul(V)\bigr)=f_\EP(1)\) and \(\EP_\chi(V,W) = \tr_W(f_\EP^V)\)
  for all admissible \(W\in\Mod_\chi(G)\).
\end{proposition}

See also \cite{Schneider-Stuhler}*{Proposition III.4.22} and
\cite{Schneider-Stuhler}*{Proposition III.4.1}.

\begin{proof}
  The finite projective resolution
  \(C_\bullet\bigl(\BT(\Gss),\gamma_e(V)\bigr)\) is explicitly built out of
  compactly induced representations.  Hence \(\Eul(V)\) belongs to the range
  of the assembly map.  Lemma~\ref{lem:truniv_cInd} yields
  \(\dim\bigl(\Eul(V)\bigr)\in\alpha\Z\) and allows us to compute \(\tr_\univ
  \Eul(V)\).  Inspection shows that this is exactly \([f_\EP^V]\).
\end{proof}

Proposition~\ref{pro:fEP} yields \(\dim V=\dim \Eul(V)\in\alpha\Z\) if~\(V\)
is irreducible supercuspidal.  This rationality result is due to Marie-France
Vignéras (\cite{Vigneras:Dimension}).

\subsection{Euler characteristics and formal dimensions for square-integrable
  representations}
\label{sec:dim_si}

\begin{theorem}
  \label{the:chi_si}
  Let \(V\in\Mod_\chi(G)\) be irreducible and square-integrable.  Let
  \[
  \iota\colon \K_0\bigl(\Hecke_\chi(G)\bigr)\to\K_0\bigl(\Sch_\chi(G)\bigr)
  \]
  be induced by the embedding \(\Hecke_\chi(G)\to\Sch_\chi(G)\).  Then
  \(\iota\bigl(\Eul(V)\bigr) = [V]\).  Hence \([V]\) lies in the range of the
  assembly map
  \[
  \varinjlim_{\mathsf{Sub}\bigl(G;C(G)\bigr)} \Rep_\chi(U)
  \to \K_0\bigl(\Hecke_\chi(G)\bigr)
  \to \K_0\bigl(\Sch_\chi(G)\bigr)
  \]
  and \(f_\EP^V(1) = \dim\bigl(\Eul(V)\bigr) = \dim^\Sch(V)\).  This number
  belongs to \(\alpha\cdot\N_{\ge1}\) with~\(\alpha\) as in
  Lemma~\ref{lem:truniv_cInd}.
\end{theorem}

\begin{proof}
  Choose \(e\in\N\) large enough such that
  \(C_\bullet\bigl(\BT(\Gss),\gamma_e(V)\bigr)\) is a projective
  \(\Hecke_\chi(G)\)\brd{}module resolution of~\(V\).  Then
  \[
  \iota \Eul(V) =
  \sum_{n=0}^\infty (-1)^n [\Sch_\chi(G)\hot_{\Hecke_\chi(G)}
  C_n\bigl(\BT(\Gss),\gamma_e(V)\bigr)].
  \]
  Since \(\Sch_\chi(G)\) is isocohomological
  (Theorem~\ref{the:Sch_chi_isocoh}), \(\Sch_\chi(G)\Lhot_{\Hecke_\chi(G)}
  V\cong V\).  Therefore, \(\Sch_\chi(G) \hot_{\Hecke_\chi(G)}
  C_\bullet\bigl(\BT(\Gss),\gamma_e(V)\bigr)\) is still a projective
  \(\Sch_\chi(G)\)\brd{}module resolution of~\(V\).  Since~\(V\) is projective
  as well (Proposition~\ref{pro:si_projective}), this resolution splits by
  bounded \(\Sch_\chi(G)\)\brd{}module homomorphisms.  This implies
  \([V]=\iota\Eul(V)\).  The remaining assertions now follow from
  Proposition~\ref{pro:fEP} and \(\dim^\Sch(V)>0\).
\end{proof}

\begin{theorem}
  \label{the:finiteness_si}
  Let \(U\in\COM\bigl(G;C(G)\bigr)\) and let \(W\in\Mod_\chi(U)\) be
  finite-dimensional.  Then there are at most \(\dim(W)\cdot\size(U)^{-1}\)
  different irreducible square-integrable representations whose restriction
  to~\(U\) contains the representation~\(W\).
\end{theorem}

\begin{proof}
  Let \(V_1,\dotsc,V_N\) be pairwise non-isomorphic irreducible
  square-integrable representations whose restriction to~\(U\) contains~\(W\).
  Let \(X\defeq \Sch_\chi(G)\hot_{\Hecke_\chi(U)} W\).  There are natural
  adjoint associativity isomorphisms
  \[
  \Hom_{\Hecke_\chi(U)}(W,V_j)\cong\Hom_{\Sch_\chi(G)}(X,V_j)
  \]
  for all~\(j\) (see~\cite{Meyer:Embed_derived}).  Thus we get non-zero maps
  \(X\to V_j\).  They are surjective and admit bounded linear sections because
  the representations~\(V_j\) are irreducible and carry the fine bornology;
  since the representations~\(V_j\) are projective
  (Proposition~\ref{pro:si_projective}), they even admit
  \(G\)\nbd{}equivariant bounded linear sections.  Thus \(V_1,\dotsc,V_N\) are
  direct summands of~\(X\).  Since they are not isomorphic, \(\bigoplus V_j\)
  is a direct summand of~\(X\) as well.  Therefore, \(\sum_{j=1}^N
  \dim^\Sch(V_j) \le \dim^\Sch X\).  Lemma~\ref{lem:truniv_cInd} and
  Theorem~\ref{the:chi_si} yield \(\dim^\Sch X = \alpha \dim(W)\size(U)^{-1}\)
  and \(\dim^\Sch(V_j)\ge\alpha\) for all~\(j\).  Hence \(N\le
  \dim(W)\size(U)^{-1}\).
\end{proof}

An irreducible square-integrable representation that is not supercuspidal is a
subquotient of a representation that we get by Jacquet induction from a proper
Levi subgroup.  It is desirable in this situation to compute the formal
dimension (and other invariants) of~\(V\) from its cuspidal data.  This gives
rise to some rather intricate computations; these are carried out
in~\cite{Aubert-Plymen:Plancherel} for representations of \(\mathrm{Gl}_m(D)\)
for a division algebra~\(D\).

\subsection{Some vanishing results}
\label{sec:vanishing}

\begin{theorem}
  \label{the:Ext_vanishing}
  Let \(V,W\in\Mod_\chi(G)\) be irreducible and tempered.  If \(V\) or~\(W\)
  is square-integrable, then \(\Ext^n_G(V,W)=0\) for all \(n\ge1\) and
  \[
  \EP_\chi(V,W)=
  \begin{cases}
    1 & \text{if \(V\cong W\);}
    \\
    0 & \text{otherwise.}
  \end{cases}
  \]
\end{theorem}

\begin{proof}
  Since \(V\) and~\(W\) are tempered and \(\Sch_\chi(G)\) is isocohomological,
  we have
  \[
  \Ext^n_{\Hecke_\chi(G)}(V,W) \cong \Ext^n_{\Sch_\chi(G)}(V,W).
  \]
  The latter vanishes for \(n\ge1\) by Proposition~\ref{pro:si_projective}.
  For \(n=0\) we are dealing with \(\Hom_G(V,W)\), which is computed by
  Schur's Lemma.
\end{proof}

\begin{theorem}
  \label{the:vanishing_dimension}
  If~\(V\) is irreducible and tempered but not square-integrable, then
  \(f_\EP^V(1) = \dim \Eul(V)=0\).
\end{theorem}

\begin{proof}
  This follows from the abstract Plancherel Theorem and
  Theorem~\ref{the:Ext_vanishing} as in the proof of
  \cite{Schneider-Stuhler}*{Corollary III.4.7}.  We merely outline the proof.
  The abstract Plancherel theorem applied to the type~I \(C^*\)\brd{}algebra
  \(C^*_{\mathrm{red},\chi}(G)\) yields that \(f_\EP^V(1)\) is the integral of
  its Fourier transform \(W\mapsto \tr_W(f_\EP^V)\) with respect to some
  measure~\(\mu\), which is called the Plancherel measure.  Here~\(W\) runs
  through the tempered irreducible representations in \(\Mod_\chi(G)\).
  Proposition~\ref{pro:fEP} asserts that \(\tr_W(f_\EP^V)=\EP_\chi(V,W)\).
  
  We have \(\Ext^n_{\Hecke_\chi(G)}(V,W)=0\) for all \(n\in\N\) and hence
  \(\EP_\chi(V,W)=0\) unless \(V\) and~\(W\) have the same infinitesimal
  character.  Since the infinitesimal character is finite-to-one, the support
  of the function \(W\mapsto \EP_\chi(V,W)\) is finite.  Hence only atoms of
  the Plancherel measure~\(\mu\) contribute to the integral
  \[
  f_\EP^V(1) = \int \EP_\chi(V,W) \,d\mu(W).
  \]
  These atoms are exactly the square-integrable representations.  Now
  Theorem~\ref{the:Ext_vanishing} yields \(f_\EP^V(1) = 0\) unless~\(V\) is
  square-integrable.  In addition, this computation shows that \(f_\EP^V(1) =
  \dim^\Sch(V)\) if~\(V\) is square-integrable (compare
  Theorem~\ref{the:chi_si}).
\end{proof}


\begin{bibdiv}
\begin{biblist}

\bib{Aubert-Plymen:Plancherel}{article}{
    author={Aubert, Anne-Marie},
    author={Plymen, Roger},
     title={Plancherel measure for \(\mathrm{GL}(n,F)\) and
       \(\mathrm{GL}(m,D)\): explicit formulas and Bernstein decomposition},
   journal={J. Number Theory},
    volume={112},
      date={2005},
    number={1},
     pages={26\ndash 66},
      issn={0022-314X},
    review={MR2131140},
}

\bib{Bass:Euler}{article}{
    author={Bass, Hyman},
     title={Euler characteristics and characters of discrete groups},
   journal={Invent. Math.},
    volume={35},
      date={1976},
     pages={155\ndash 196},
    review={MR0432781 (55 \#5764)},
}

\bib{Bernstein:centre}{article}{
    author={Bernstein, J. N.},
     title={Le ``centre'' de Bernstein},
      book={
       title={Representations of reductive groups over a local field},
    language={French},
   publisher={Hermann},
       place={Paris},
        date={1984},
      series={Travaux en Cours},
        },
     pages={1\ndash 32},
      note={Edited by P. Deligne},
    review={MR771671 (86e:22028)},
}

\bib{Bridson-Haefliger}{book}{
    author={Bridson, Martin R.},
    author={Haefliger, André},
     title={Metric spaces of non-positive curvature},
    series={Grundlehren der Mathematischen Wissenschaften},
    volume={319},
 publisher={Springer-Verlag},
     place={Berlin},
      date={1999},
     pages={xxii+643},
      isbn={3-540-64324-9},
    review={MR1744486 (2000k:53038)},
}

\bib{Bruhat-Tits:Groupes_reductifs}{article}{
    author={Bruhat, F.},
    author={Tits, J.},
     title={Groupes réductifs sur un corps local},
  language={French},
   journal={Inst. Hautes Études Sci. Publ. Math.},
    number={41},
      date={1972},
     pages={5\ndash 251},
    review={MR0327923 (48 \#6265)},
}

\bib{Emerson-Meyer:Dualizing}{article}{
    author={Emerson, Heath},
    author={Meyer, Ralf},
     title={Dualizing the coarse assembly map},
      date={2004},
    status={to appear},
    eprint={http://arxiv.org/math.OA/0401227},
   journal={J. Inst. Math. Jussieu},
}

\bib{Grothendieck:Produits}{book}{
    author={Grothendieck, Alexander},
     title={Produits tensoriels topologiques et espaces nucléaires},
  language={French},
    series={Mem. Amer. Math. Soc.},
    volume={16},
      date={1955},
     pages={140},
    review={MR0075539 (17,763c)},
}

\bib{Hogbe-Nlend:Completions}{article}{
    author={Hogbe-Nlend, Henri},
     title={Complétion, tenseurs et nucléarité en bornologie},
  language={French},
   journal={J. Math. Pures Appl. (9)},
    volume={49},
      date={1970},
     pages={193\ndash 288},
    review={MR0279557 (43 \#5279)},
}

\bib{Hogbe-Nlend:Bornologies}{book}{
    author={Hogbe-Nlend, Henri},
     title={Bornologies and functional analysis},
    series={North-Holland Mathematics Studies},
    volume={26},
 publisher={North-Holland Publishing Co.},
     place={Amsterdam},
      date={1977},
     pages={xii+144},
      isbn={0-7204-0712-5},
    review={MR0500064 (58 \#17774)},
}

\bib{Humphreys:Coxeter_groups}{book}{
    author={Humphreys, James E.},
     title={Reflection groups and Coxeter groups},
    series={Cambridge Studies in Advanced Mathematics},
    volume={29},
 publisher={Cambridge University Press},
     place={Cambridge},
      date={1990},
     pages={xii+204},
      isbn={0-521-37510-X},
    review={MR1066460 (92h:20002)},
}

\bib{Keller:Handbook}{incollection}{
    author={Keller, Bernhard},
     title={Derived categories and their uses},
      date={1996},
      book={
       title={Handbook of algebra, vol.\ 1},
   publisher={North-Holland},
     address={Amsterdam},
       },
     pages={671\ndash 701},
    review={MR1421815 (98h:18013)},
}

\bib{Lueck-Reich:Survey_assembly}{article}{
    author={Lück, Wolfgang},
    author={Reich, Holger},
     title={The Baum-Connes and the Farrell-Jones Conjectures in K- and
            L\nobreakdash-theory},
   journal={Preprintreihe SFB 478},
      note={Universität Münster},
    volume={324},
      date={2004},
}

\bib{Meyer:Analytic}{thesis}{
    author={Meyer, Ralf},
     title={Analytic cyclic cohomology},
      type={Ph.D. Thesis},
institution={Westfälische Wilhelms-Universität Münster},
      date={1999},
    eprint={http://arxiv.org/math.KT/9906205},
}

\bib{Meyer:Born_Top}{article}{
    author={Meyer, Ralf},
     title={Bornological versus topological analysis in metrizable spaces},
      book={
       title={Banach algebras and their applications},
      series={Contemporary Mathematics},
      volume={363},
      editor={Anthony To-Ming Lau},
      editor={Volker Runde},
   publisher={American Mathematical Society},
       place={Providence, RI},
        date={2004},
      },
     pages={249\ndash 278},
    review={MR2097966},
}

\bib{Meyer:Smooth}{article}{
    author={Meyer, Ralf},
     title={Smooth group representations on bornological vector spaces},
  language={English, with English and French summaries},
   journal={Bull. Sci. Math.},
    volume={128},
      date={2004},
    number={2},
     pages={127\ndash 166},
      issn={0007-4497},
    review={MR2039113 (2005c:22013)},
}

\bib{Meyer:Poly_comb}{article}{
    author={Meyer, Ralf},
     title={Combable groups have group cohomology of polynomial growth},
      date={2004},
    status={to appear in Q. J. Math.},
    eprint={http://arxiv.org/math.KT/0410597},
}

\bib{Meyer:Embed_derived}{article}{
    author={Meyer, Ralf},
     title={Embeddings of derived categories of bornological modules},
      date={2004},
    status={eprint},
    eprint={http://arxiv.org/math.FA/0410596},
}

\bib{Pirkovskii:Stably_flat}{article}{
    author={Pirkovskii, A. Yu.},
     title={Stably flat completions of universal enveloping algebras},
      date={2003},
    status={eprint},
    eprint={http://arxiv.org/math.FA/0311492},
}

\bib{Rosenberg:Algebraic_K}{book}{
    author={Rosenberg, Jonathan},
     title={Algebraic K-theory and its applications},
    series={Graduate Texts in Mathematics},
    volume={147},
 publisher={Springer-Verlag},
     place={New York},
      date={1994},
     pages={x+392},
      isbn={0-387-94248-3},
    review={MR1282290 (95e:19001)},
}

\bib{Schneider-Stuhler}{article}{
    author={Schneider, Peter},
    author={Stuhler, Ulrich},
     title={Representation theory and sheaves on the Bruhat-Tits building},
   journal={Inst. Hautes Études Sci. Publ. Math.},
    number={85},
      date={1997},
     pages={97\ndash 191},
      issn={0073-8301},
    review={MR1471867 (98m:22023)},
}

\bib{Schneider-Zink}{article}{
    author={Schneider, P.},
    author={Zink, E.-W.},
     title={\(K\)-types for the tempered components of a \(p\)-adic general
       linear group},
      note={With an appendix by P. Schneider and U. Stuhler},
   journal={J. Reine Angew. Math.},
    volume={517},
      date={1999},
     pages={161\ndash 208},
      issn={0075-4102},
    review={MR1728541 (2001f:22029)},
}

\bib{Silberger:Reductive_padic}{book}{
    author={Silberger, Allan J.},
     title={Introduction to harmonic analysis on reductive \(p\)-adic groups},
    series={Mathematical Notes},
    volume={23},
      note={Based on lectures by Harish-Chandra at the Institute for
            Advanced Study, 1971--1973},
 publisher={Princeton University Press},
     place={Princeton, N.J.},
      date={1979},
     pages={iv+371},
      isbn={0-691-08246-4},
    review={MR544991 (81m:22025)},
}

\bib{Tits:Corvallis}{article}{
    author={Tits, J.},
     title={Reductive groups over local fields},
      book={
       title={Automorphic forms, representations and
           \(L\)\nobreakdash-functions},
      series={Proc. Sympos. Pure Math., XXXIII,Part 1},
   publisher={Amer. Math. Soc.},
       place={Providence, R.I.},
        date={1979},
        },
conference={
       title={Sympos. Pure Math.},
       place={Oregon State Univ., Corvallis, Ore.},
        date={1977},
        },
     pages={29\ndash 69},
    review={MR546588 (80h:20064)},
}

\bib{Treves:Kernels}{book}{
    author={Trèves, François},
     title={Topological vector spaces, distributions and kernels},
 publisher={Academic Press},
     place={New York},
      date={1967},
     pages={xvi+624},
    review={MR0225131 (37 \#726)},
}

\bib{Vigneras:Dimension}{article}{
    author={Vignéras, Marie-France},
     title={On formal dimensions for reductive \(p\)-adic groups},
      book={
       title={Festschrift in honor of I. I. Piatetski-Shapiro on the
         occasion of his sixtieth birthday, Part I (Ramat Aviv, 1989)},
      series={Israel Math. Conf. Proc.},
      volume={2},
 publisher={Weizmann},
     place={Jerusalem},
      date={1990},
      },
     pages={225\ndash 266},
    review={MR1159104 (93c:22034)},
}

\bib{Waldspurger:Plancherel_formula}{article}{
    author={Waldspurger, J.-L.},
     title={La formule de Plancherel pour les groupes \(p\)-adiques (d'après
            Harish-Chandra)},
  language={French},
   journal={J. Inst. Math. Jussieu},
    volume={2},
      date={2003},
    number={2},
     pages={235\ndash 333},
      issn={1474-7480},
    review={MR1989693 (2004d:22009)},
}

\end{biblist}
\end{bibdiv}
\end{document}